\documentclass[3p,english]{elsarticle}

\usepackage[utf8]{inputenc}
\usepackage[T1]{fontenc}

\usepackage{babel}
\usepackage{tikz}
\usepackage{pgfplots}

\usepackage{graphicx}
\usepackage{amsmath,amsfonts,amssymb}

\usepackage{xspace}
\usepackage{fourier}
\usepackage{geometry}

\usepackage{titlesec}

  \usepackage[colorlinks=true,linkcolor=black,citecolor=blue]{hyperref}

\DeclareMathOperator{\dv}{div}

\newcommand{\R}{\mathbb{R}}
\newcommand{\N}{\mathbb{N}}

\newcommand{\x}{\times }

\newcommand{\eps}{\varepsilon}
\newcommand{\disp}{\displaystyle}

\newcommand{\ov}{\overline}

\newcommand{\deriv}[2]{ {\frac{\partial #1}{\partial #2} }}
\newcommand{\hsoil}{ h_{\mathrm{soil}} }
\newcommand{\hbot}{ h_{\mathrm{bot}} }
\newcommand{\hsat}{ h_{\mathrm{sat}} }
\newcommand{\hsatdd}{ \hsat^{ 2d} }
\newcommand{\hsatmin}{ \hsat^{ a} }
\newcommand{\hsatmax}{ \hsat^{ b} }
\newcommand{\hsatint}{ \hsat^{ c} }

\newcommand{\hsoilb}{ \ov h_{\mathrm{soil}} }
\newcommand{\hbotb}{ \ov h_{\mathrm{bot}} }
\newcommand{\hmaxb}{ \ov h_{\mathrm{max}} }

\newcommand{\Gammat}{ \Gamma_{\mathrm{soil}} }
\newcommand{\Gammas}{\Gammat}
\newcommand{\Gammab}{ \Gamma_{\mathrm{bot}} }
\newcommand{\Gammav}{ \Gamma_{\mathrm{ver}} }

\newcommand{\Gammatb}{\ov \Gamma_{\mathrm{soil}} }
\newcommand{\Gammabb}{ \ov\Gamma_{\mathrm{bot}} }
\newcommand{\Gammavb}{ \ov\Gamma_{\mathrm{ver}} }

\newcommand{\Pinit}{ P_{\mathrm{init}} }
\newcommand{\Hinit}{ H_{\mathrm{init}} }

\newcommand{\hmax}{ h_{\mathrm{max}} }

\newcommand{\fluxsoil}{ F }
\newcommand{\ha}{ h}

\newcommand{\hab}{ \ov h}

\newcommand{\Omplus}{ {\Omega_{\ha}^+ }}
\newcommand{\Ommoins}{ {\Omega_{\ha}^- }}
\newcommand{\Ommoinsat}{ {\Omega_{\hsat}^- }}

\newcommand{\Omplusb}{ {\Omega_{\hab}^+ }}
\newcommand{\Ommoinsb}{ {\Omega_{\hab}^- }}

\newcommand{\interface}{ {\Gamma_{\ha}} }
\newcommand{\interfacep}{ {\Gamma_{\ha}^+} }
\newcommand{\interfaceb}{ {\Gamma_{\hab}} }

\newcommand{\Mr}{ {M_{\ha}^+ }}
\newcommand{\Md}{ {M_{\ha}^- }}

\newcommand{\Pplus}{P }

\newcommand{\vplus}{ {u }}

\newcommand{\Mtot}{ M_{\mathrm{tot}} }
\newcommand{\domt}{ ]0,T[ }
\newcommand{\domtb}{ ]0,\ov T[ }

\newcommand{\Pbub}{P_{s}}
\newcommand{\Pbubb}{\ov P_{s}}
\newcommand{\refmain}{\eqref{coupled_transition}--\eqref{coupled_velocity}\xspace}
\newcommand{\refrescaled}{\eqref{dimensionless_transition}--\eqref{dimensionless_interface}\xspace}

\renewcommand{\and}{\quad\text{and}\quad }

\newcommand{\qand}{\qquad\text{and}\qquad }

\newtheorem{theorem}{Theorem}[section]
\newtheorem{proposition}[theorem]{Proposition}

\newdefinition{remark}{Remark}

\newproof{pf}{Proof}

\numberwithin{equation}{section}

\pgfplotsset{compat=1.10}

\author[C. Bourel]{Christophe Bourel\corref{cor1}}
\author[C. Choquet]{Catherine Choquet}
\author[C. Bourel]{Carole Rosier}
\author[C. Bourel]{Munkhgerel Tsegmid}

\cortext[cor1]{Corresponding author: christophe.bourel@univ-littoral.fr}

\title{Modelling of shallow aquifers in interaction with overland water}

\titleformat{\paragraph}[runin]{\bfseries}{}{}{}
\titleformat{\section}[block]{\bfseries\scshape}{\thesection.}{0.4em}{}
\titleformat{\subsection}[block]{\bfseries}{\thesubsection.}{0.4em}{}

\address[C. Bourel]{Univ. Littoral C\^ote d'Opale, EA 2797 - LMPA, F- 62228 Calais, France}
\address[C. Choquet]{La Rochelle Universit\'e,  MIA, EA 3165,  F-17031 La Rochelle, France}

\begin{document}

\begin{abstract}
In this work, we present a class of new efficient models  for  water flow in shallow unconfined aquifers, giving an alternative to the classical but less tractable 3D-Richards model.
Its derivation is guided by two ambitions: any new model should be low cost in computational time and should still give relevant results at every time scale.
We thus keep track of  two types of flow occurring in such a context and which are dominant when the {\it ratio} thickness over longitudinal length is small: 
the first one is dominant in a small time scale and is described by a 
vertical 1D-Richards problem; the second one corresponds to a large time scale, when the evolution of the hydraulic head  turns to become independent of  the vertical variable. 
These two types of flow are appropriately modelled by, respectively, a one-dimensional and a two-dimensional system of PDEs boundary value problems. 
They are coupled along an artificial level below which the Dupuit hypothesis holds true ({\it i.e.} the vertical flow is instantaneous) in a way ensuring that the global model is mass conservative.
Tuning the artificial level, which even can depend on an unknown of the problem, we browse  the new class of models. 
We prove using asymptotic expansions that the 3D-Richards problem and each model of the class behaves the same at every considered time scale (short, intermediate and large) in thin aquifers. 
The results are illustrated by numerical simulations, showing especially that the new models results fit well with the ones obtained with the  original 3d-Richards problem even in non-thin aquifers.
\end{abstract}


\maketitle


\section{Introduction}
Contamination of soil and groundwater is a major concern that affects all populated areas. 
Many works are thus  developed for studying the vulnerability
of aquifers with regard to agricultural, industrial,  or sewage pollutions. 
There is an abundant literature on each of the involved processes (geological, physical, chemical...), so that we can consider that the corresponding model is already available.
Nevertheless there is a so wide variety of processes (chemical, hydrogeological, anthropic)
acting in a so wide range of temporal and geometrical length scales that the assembly
of the corresponding model bricks, if considered like toolboxes of a software, is, at best,
computationally expensive. 

In this multi-scale context, a particularly interesting issue is a proper and tractable model for the exchanges between the overland and the underground waters. 
Indeed, the challenge consists in capturing very different physical phenomena, the fast and essentially vertical leakage coming from the surface through an unsaturated soil and the slow and essentially horizontal 
displacement in the saturated part of the aquifer, that are classically modelled by mathematical systems with very different structures. 
The question is all the more important that an accurate study of the interaction between  the water table and the overland water is essential for many concerns, concerns that disallow the use of classical time upscaling 
processes. 
It is in particular crucial for studying the transport of chemical components in the aquifer. Indeed, it turns out that many chemical reactions occur in the first meters of the subsoil, where oxygen is still very present. 
As a byproduct, the chemical species that reach the water table are not necessarily the same than  those that have left the surface, and there is a large range of kinetics reaction times to handle with. 
There is actually no scale separation.

In the present paper, we focus on the hydrogeological question. 
We thus consider the displacement of a wetting phase (water) in the presence of a non-wetting fluid (air) in a porous medium.
Assuming that the air present in the unsaturated zone has infinite mobility allows to use a model for immiscible fluid flow simplified by the Richards hypothesis.
The saturation is thus considered as a monotone function depending of the pressure head and the so-called Richards model consists in a nonlinear  three-dimensional equation of degenerate parabolic type.
All the existing simplified models for the fluid displacement in aquifers are motivated by the characteristics of the flow in their saturated part.
A form of stratification enables the definition of interfaces and  the slowness of the natural dynamics ensures that these interfaces have a smooth and stable behaviour. Moreover the flows are essentially orthogonal to the walls (Dupuit's hypothesis). 
These points allow the vertical integration of the Richards equation in the saturated area and lead to the use of a family of 2D models developed since the 60's (see e.g. the works of Jacob Bear,  \cite{Bear1972,Bear1987}). 
A main weakness of the approach by vertical integration lies in its justification. 
It is only valuable for very precise length and time scales, the time
scale in particular being completely different of the typical durations of chemical reactions (see once again \cite{Bear1972} for empirical  and qualitative arguments, see \cite{Pan}  for asymptotic computations). 
However, such 2D models are widely used, even out of their validity range and even if it turns out to be  especially difficult to properly couple them with the flow in the unsaturated part of the underground.
Only numerical attempts were done in this direction.
We mention \cite {Kong} where the integrated model is directly coupled with a surface model (see also the references therein). 
The unsaturated area of the aquifer is taken into account  in \cite{pikul1974numerical} using a $1D$-Richards equation coupled with a simplified model in the saturated part.
However, the study is purely numerical and the model is not mathematically justified. In \cite{abbott1986introduction}, the latter 
kind of model is integrated into a computational code called "SHE" (for "European Hydrological
System" and later became SHETRAN) in the case where the water table remains away
from the ground level. 
See also \cite{Yak}, \cite{Paulus}.

To the best of our knowledge, there is no mathematical justification for any
 "Dupuit-Richards" model specifying the hypotheses as well as the scales that allow
its derivation from a more complete model (such as the $3D$-Richards one). 

Notice finally that the coupling of the surface and underground flows turns out to be more tractable when handling with a Richards equation (see e.g. \cite{Ern}  
or \cite{Ber} and \cite{Bern} where  the surface behavior is reduced to a Signorini boundary condition).

\bigskip

The goal of this work is to provide a simple model exploiting the low thickness of a confined or unconfined aquifer.
In summary it consists in coupling purely vertical models (describing the flow at a small time scale) with a horizontal model (describing the flow at a long time scale). 
Clearly, given its construction, the model is simpler to manipulate numerically since the original $3D$ problem is replaced by the coupling 
of a $2D$ problem with several independent $1D$-problems (which can be solved in parallel).
Significant time savings are expected in the numerical processing.

This work could be viewed as another attempt using the numerically pragmatical methodology of \cite{abbott1986introduction} and leading to a ``Dupuit-Richards'' model.
Yet, our approach is quite different.
First, we actually derive a class of models, each of them being characterised by the definition of some virtual interface which does not necessarily  coincide with the water table (especially when trying to optimize the error).
It follows that a model of this class does not necessarily contains a Dupuit component. 
The position of the virtual interface may even be an unknown of our model. 
Next, we aim at describing the flow in a large range of time scales, and, more precisely without any assumption of scale separation.
The idea consists in always capturing  both the fast and slow components of the flow  given by Richards $3D$ equations, whatever the time scale.  
Their coupling is done through flux terms ensuring that the model is  mass conservative (and thus avoiding the criticism done in \cite{Vachaud}). 
Finally, the large validity range of the new class of models is justified by an asymptotic study.
But, as already mentioned, no time scale separation is assumed in the present paper so that we adopt a new methodology for the asymptotic arguments.
Let $\eps >0$ describe the ratio of the aquifer's deepness over its characteristic horizontal length. 
Assume that $\eps$ is small. 
The usual approach would consist in choosing a reference time for the study, introducing an asymptotic expansion of the solution of the 3D-Richards system and using the scale separation for identifying the equations governing the main order terms of this ansatz.
This is the classical process for deriving an effective model. 
Here the asymptotic analysis is not used for deriving an effective model for a given reference time.
Rather, it is used for proving that each model of our new class and the 3D-Richards equation are  associated with the same effective problem for any time scale.
Basically:  
\begin{enumerate}
 \item At  short times, the horizontal flow is very small and the vertical one satisfies a  $1D$-Richards problem.
 \item At non-short times, the vertical flow appears  instantaneous. The corresponding pressure profile satisfies the stationary $1D$-Richards problem. Then the hydraulic head $H$ does not depend on the vertical variable $z$. This corresponds to the so-called Dupuit hypothesis.
  \item At large times, the horizontal flux is non-zero. It is ruled by a $2D$-horizontal diffusion equation where the conductivity  is the vertical average of the permeability tensor on the \emph{whole} depth of the aquifer.
\end{enumerate}

\bigskip

The paper is organised as follows: In Section \ref{sec_intro},  we describe the geometry of the problem, the physical parameters and unknowns.
The classical 3D-Richards model is recalled. 
The main result and numerical simulations are given in Section  \ref{sec_main}. Namely, we present the systems coupling the vertical  and the horizontal flows and we comment 
on the model. Finally, the formal asymptotic analysis of our models and of the 3D-Richards model are performed and compared  in Section \ref{formal_asymptotic}.


\section{Description of the problem}\label{sec_intro}

This section is devoted to the description of the domain of study, of the physical parameters and of the unknowns which are chosen for characterising the flow through the Richards model.


\subsection{Geometry}


The aquifer corresponds to a cylindrical domain $\Omega\subset\R^3$.
For the sake of the simplicity, we assume vertical walls. 
The projection of $\Omega$ on any horizontal  plane is an open domain $\Omega_x\subset\R^2$ with boundary $\partial\Omega_x$. 
The lower and upper  bases of $\Omega$ are respectively the graphs of real-valued functions $\hbot$ and $\hsoil$ such  that
\begin{equation}\label{soilbot}
\hsoil(x)>\hbot(x)\ ,\qquad \forall x\in\Omega_x  .
\end{equation} 
In summary the domain is given by:
\begin{equation}\label{omega}
\Omega=\big\{(x,z)\in\Omega_x\times\R\quad |\quad z\in\big]\hbot(x) ,\hsoil(x)\big[    \big\} .
\end{equation}
We split the boundary $\partial\Omega$ of $\Omega$  in three parts (bottom, top and vertical)
\begin{gather*}
\partial\Omega = \Gammab\sqcup\Gammas\sqcup\Gammav\, , 
\\
\Gammab:=\big\{(x,z)\in\Omega\ |\ z=\hbot(x) \big\}\ , \quad  \Gammas:=\big\{(x,z)\in\Omega\ |\ z=\hsoil(x)\big\}\ , \quad  \Gammav:=\big\{(x,z)\in\Omega\ |\ x\in\partial\Omega_x\big\}.
\end{gather*}
In the present paper, as already mentioned,  we derive a class of models that are characterised by the position $h$ of some virtual interface in the reservoir.
For our construction, this function has to take its values  in the semi-open interval $[\hbot,\hsoil)$. 
For numerical implementation, an easy recipe consists in replacing the condition $h < \hsoil$ by $h \le \hsoil-\delta$ where $\delta $ is an arbitrary  small positive real number. 
We thus introduce the auxiliary function $\hmax$ defined by
\begin{equation}\label{hamax}
\hmax=\hsoil-\delta , \quad 0<\delta \ll 1.
\end{equation}

\begin{figure}\label{fig.geometry}
	\begin{center}
\includegraphics[scale=0.4]{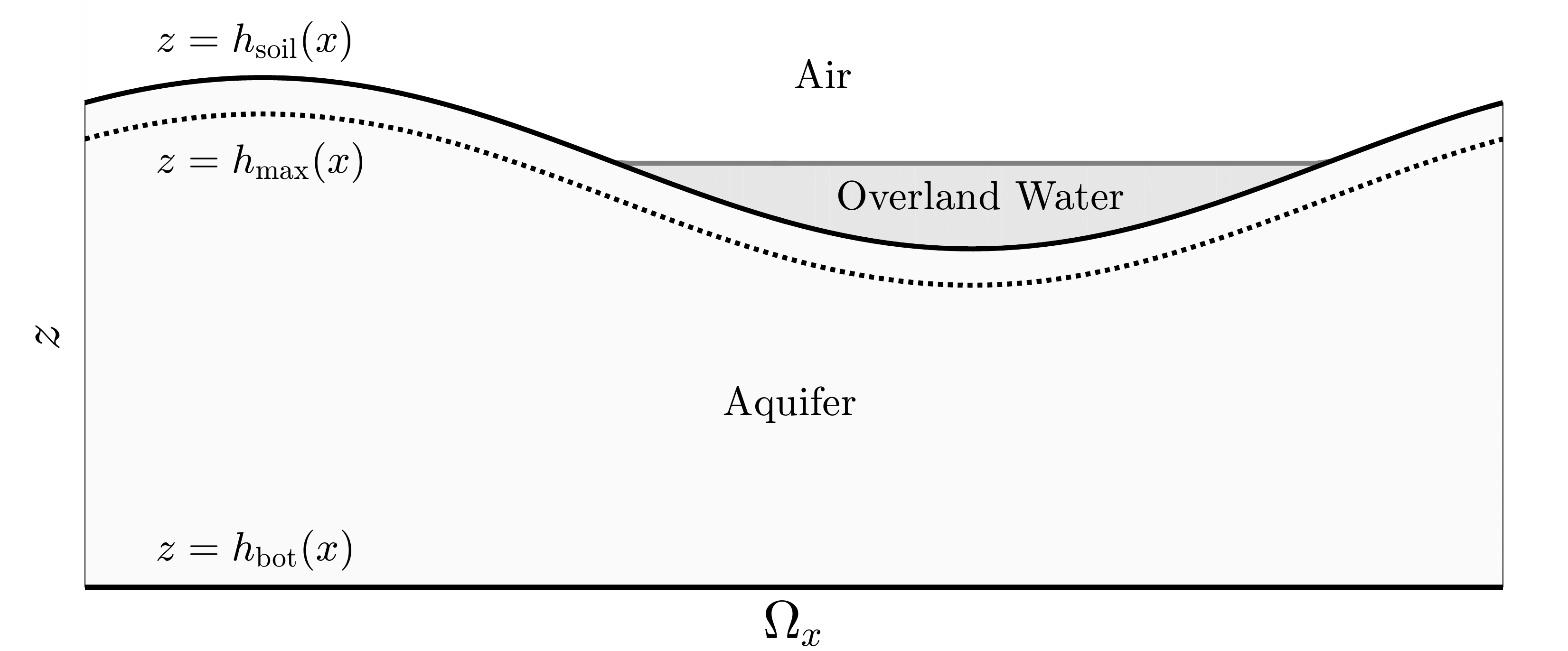}
		\caption{ Bidimensional representation of the cylindrical geometry of the problem: $\Omega_x\subset\R$ is an interval.   \label{pres_geo}}
	\end{center}
\end{figure}


\subsection{Three-dimensional Richards equation}\label{sec.model3d}

\begin{figure}
	\begin{minipage}{0.49\linewidth}
		\begin{tikzpicture}
		\draw[->] (-4.5,0) -- (1.7,0) node[below left] {$P$};
		\draw[->] (0,-0.5) -- (0,1.5) node[left] {$s(P)$};
		\draw[scale=1,domain=-4.5:-1.5,smooth,variable=\x,blue,line width=0.9pt] plot ({\x},{(-1.5/\x)^3});
		\draw[scale=1,domain=-1.5:1.5,smooth,variable=\x,blue,line width=0.9pt] plot ({\x},{1});
		\draw[-,dashed] (-1.5,0)node[below] {$P_s$} -- (-1.5,1) ;
		\end{tikzpicture}
	\end{minipage}
	\hfill
	\begin{minipage}{0.49\linewidth}
		\begin{tikzpicture}
		\draw[->] (-4.5,0) -- (1.7,0) node[below left] {$P$};
		\draw[->] (0,-0.5) -- (0,1.5) node[left] {$k_r(P)$};
		\draw[scale=1,domain=-4.5:-1.5,smooth,variable=\x,blue,line width=0.9pt] plot ({\x},{(-1.5/\x)^11});
		\draw[scale=1,domain=-1.5:1.5,smooth,variable=\x,blue,line width=0.9pt] plot ({\x},{1});
		\draw[-,dashed] (-1.5,0)node[below] {$P_s$} -- (-1.5,1) ;
		\end{tikzpicture}
		
	\end{minipage}
	
	\caption{Saturation and relative permeability in terms of the pressure: the Brooks and Corey model.  \label{saturation_profile_figure}}
\end{figure}
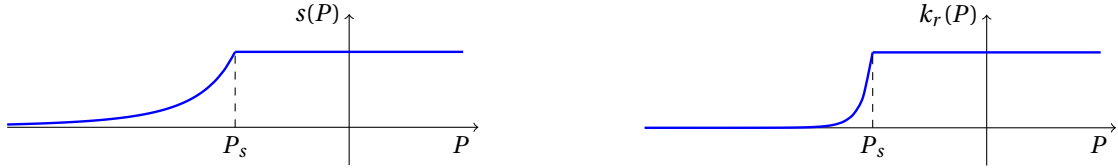

We aim at deriving alternatives to the Richards equation. 
Let us briefly describe this classical model.  
In this paper we limit our study to a one-phase incompressible fluid which accordingly admits
a constant density $\rho\in \R^{*}_+$.
First, in multiphase systems, observations have shown that an increase of the saturation of the non-wetting phase leads to an increase of the capillary pressure.
The Richards model is moreover based on the assumption that the air pressure in the underground equals the atmospheric pressure, thus is not an unknown of the problem.
One thus assumes that the saturation and the relative conductivity of the soil are given as \textit{functions} of the fluid pressure $P$, denoted respectively by $s=s(P)$ and $k_r=k_r(P)$. 
There is a large choice of available models for~$s$ and $k_r$. The most classical  examples for an air-water system are the van Genuchten model \cite{van1980closed}, with no-explicit dependance on the bubbling pressure but with fitting parameters, and the Brooks and Corey model \cite{brooks1964hydraulic}, that we use in the simulations below:
	\begin{equation}\label{brooks}
	s(P)=\begin{cases}
	\left(\Pbub/P\right)^\lambda&\text{if }P<\Pbub\\
	1&\text{if }P\ge\Pbub\\
	\end{cases},\qquad
	k_r(P)=\begin{cases}
	\left(\Pbub/P\right)^\gamma&\text{if }P<\Pbub\\
	1&\text{if }P\ge\Pbub\\
	\end{cases},
	\end{equation}
	where $\lambda>0$, $\gamma=2+3\lambda$ and $\Pbub<0$.
Notice that our model would easily adapt to  hysteretic soil properties (\cite{Pham}, \cite{Schw}).	
Since these methods, as of today, do not permit three-dimensional calculations, we guess that our 1D-2D models are even more interesting for their implementation than the 3D-Richards model. 
The important point is that these models are such that 
\begin{equation}\label{sat}
s(P)=1 \quad\Longleftrightarrow\quad P\ge\Pbub\qand k_r(P)=1 \quad\Longleftrightarrow\quad P\ge\Pbub.
\end{equation}
In particular, the water pressure is greater than the bubbling pressure $\Pbub$ if and only if the soil is completely saturated ($\Pbub$ being a fixed real number).
The graphs of the functions $s$ and $k_r$ given by the Brooks-Corey model used below for the numerical simulations are represented in Figure \ref{saturation_profile_figure} (the parameters are given at the beginning of Subsection \ref{numeric}).

The soil transmission properties are characterised by the porosity function, $\phi=\phi(x,z)\in(0,1)$, and the  permeability tensor,  $K_0(x,z)$. 
The latter  is a $3\times3$ symmetric positive definite tensor which describes the conductivity of the \textit{saturated} soil at the position $(x,z)\in\Omega$. 
We introduce $K_{xx}\in \mathcal{M}_{22}(\R)$, $K_{zz}\in \R^*$ and $K_{xz}\in\mathcal{M}_{21}(\R)$ such that
\begin{equation}\label{K0}
K_0=\begin{pmatrix}
K_{xx}&K_{xz}\\
K_{xz}^T&K_{zz}
\end{pmatrix}.
\end{equation}

The fluid is characterised by its pressure $P$ and its velocity $v$ solving the following Richards problem:
\begin{equation}\label{richards3d}
\begin{cases}
\disp\phi \deriv{s(P)}{t}+\dv(v)=0&\text{in }\domt\x\Omega\\
\disp v = -k_r(P)\,K_0\,\Big( \frac{1}{\rho g}\nabla P+e_3 \Big)&\text{in }\domt\x\Omega\\
\alpha\, P + \beta\, v\cdot n = F&\text{on }\domt\x\Gammat\\
v\cdot n = 0&\text{on }\domt\x(\Gammab\cup\Gammav)
\end{cases}
\end{equation}
where $g$ is the gravity constant and $e_3$ is the unitary vertical vector pointing up.
The first equation describes the mass conservation of the constant density fluid in the case of an incompressible soil. The second equation is the Darcy's law associated with the nonlinear anisotropic conductivity $k_r(P)\,K_0$. 
The boundary condition $v\cdot n=0$ on $\Gammab$ corresponds to the impermeable layer at the bottom of the aquifer. 
The same is assumed on   $\Gammav$  to  simplify the presentation.
The condition at the soil level $\Gammas$ is a Robin condition associated with given $(\alpha,\beta)\in(\R_+)^2\setminus \{0,0\}$ and $F:\Gammas\to \R$.

\begin{remark}[Dominant behaviors in a shallow aquifer]\label{rem_behavior}
 In Section \ref{formal_asymptotic} we investigate the behavior of the flow described by the 3D-Richards equations in the case of a thin aquifer and for various time scales. 
Let us summarise the conclusions of this asymptotic analysis.
They might shed light on the comments about our models in the next section.
\begin{enumerate}
 \item At any time scale, the dominant flow is the one in the vertical direction (see for example \eqref{rescaled_mass} in which the horizontal diffusion term appears multiplied by the small parameter~$\eps$).
 \item In the short time scale ($T\sim 1$), the horizontal flow is very small and the vertical one solves a classical 1D-Richards problem.
 \item In non-short time scales ($T\sim \eps^{-1}$ or $T\sim \eps^{-2}$), the vertical flow appears as being instantaneous. The corresponding pressure profile satisfies  a stationary 1D-Richards problem. Then the pressure is $P=\rho\,g(H-z)$ where the hydraulic head $H$ does not depend on the vertical variable $z$. The  velocity is horizontal.  This corresponds to the so-called Dupuit hypothesis.
 \item In the long time scale ($T\sim \eps^{-2}$), the horizontal flow is non-zero and it is ruled by a 2D-horizontal diffusion equation where the conductivity  is the vertical average of the permeability tensor on the \emph{whole} depth of the aquifer, from $\hbot$ to $\hsoil$.
\end{enumerate}
\end{remark}


\section{Main result and numerical simulations }\label{sec_main}


\subsection{Models coupling vertical 1d-Richards flow and Dupuit horizontal flow}


Each of our models  splits the description of the flow into two subregions of $\Omega$ (possibly time-dependent). 
These  zones are defined by  a function $\ha=\ha(t,x)$  such that $\hbot\leq\ha<\hsoil$:
 \begin{gather}
\Ommoins(t):=  \big\{ (x,z)\in\Omega\ |\ z<\ha(x,t) \big \}\quad\text{and}\quad\Omplus(t):= \big\{ (x,z)\in\Omega\ |\ z>\ha(x,t) \big \} ,
\label{omplusmoins}
\\
\interface:=\big\{ (x,z)\in\Omega\ |\ z=\ha(x,t) \big \} .
\label{gamaha}
\end{gather}
We emphasise  that choosing the level $\ha$ corresponds to the specification of one of the models of our class.
The function $\ha$  can even be an unknown of our problem, more precisely depending of an unknown of the problem (see condition \eqref{coupled_interface} below).

On the other hand we introduce the following tensor $M_0$ which will act as an effective permeability tensor:
\begin{equation}\label{M0}
M_0=\begin{pmatrix}
                                                                        S_0&0\\0&0
                                                                       \end{pmatrix},
\qquad
S_0=K_{xx}-\frac{1}{K_{zz}}K_{xz}K_{zx} .
\end{equation}
The $2\x2$ matrix $S_0$ is the Schur complement of the block $K_{zz}$ in the tensor $K_0$.
Since $K_0$ is a symetric positive definite matrix (see just before \eqref{K0}),  the same holds  for $S_0$.
We then introduce the averaged conductivity tensor $\tilde K$ defined in $\domt\x\Omega_x$ for any function $\tilde H=\tilde H(t,x)$ by
\begin{equation}\label{averaged_conductivity}
\tilde K(\tilde H)(t,x)=\int_{\hbot(x)}^{\hsoil(x)}k_r\big(\rho\,g(\tilde H(t,x)-z)\big)\,M_0(x,z)\,dz.
\end{equation}

Finally, for the 2D part of the model, we  introduce the  notations $\nabla_x=(\partial_{x_1}, \partial_{x_2}, 0)^T$ for  the horizontal gradient and $\dv_x(v)=\nabla_x\cdot v=\partial_{x_1}v_1+\partial_{x_2}v_2$ for the horizontal divergence of $v\in \R^3$.

\paragraph{The model.}
Our coupled model consists in finding the pressure $P$, the velocity $v$ and the auxiliary unknowns  $u$, $w$, $\tilde H$  and $\ha$ such that:
\begin{itemize}
 \item In $\Omplus(t)$ the following 1D-Richards equation holds
\begin{equation}\label{coupled_transition}
\begin{cases}
\disp \phi \deriv{s(P)}{t}+\deriv{}{z}\big(u\cdot e_3\big)=0 & \text{for \ }t\in]0,T[\ ,\quad (x,z)\in\Omplus(t)\\
 \alpha\, P + \beta\, u\cdot e_3=F  &\text{for \ }(t,x)\in]0,T[\x\Gammat\\
  P\big(t,x,\ha(t,x)\big)=\rho\,g\big(\tilde H(t,x)-\ha(t,x) \big)&\text{for \ }(t,x)\in]0,T[\x\Omega_x\\
  P(0,x,z)=\Pinit(x,z)&\text{for \ } (x,z)\in\Omplus(0)\\
\end{cases}
\end{equation}
\item In  $\Ommoins(t)$ the  pressure $P$ satisfies
\begin{equation}\label{coupled_watertable}
       \disp P(t,x,z)=\rho\,g\big( \tilde H(t,x)-z \big) \qquad \text{for \ }t\in[0,T[\ ,\quad (x,z)\in\Ommoins(t)\\
\end{equation}

\item The hydraulic head solves in $\Omega_x$
\begin{equation}\label{coupled_hydraulic}
\begin{cases}
\dv_x\Big(\tilde K(\tilde H)\,\nabla_x\tilde H \Big)=(u\cdot e_3)\big|_{\interfacep}&\text{for \ }(t,x)\in]0,T[\x\Omega_x
 \\
 \disp  \tilde K(\tilde H)\,\nabla_x \tilde H\cdot n=0&\text{for \ }(t,x)\in]0,T[\x\partial\Omega_x\\
  \tilde H(0,x)=\Hinit(x)&\text{for \ } x\in\Omega_x
\end{cases}
\end{equation} 
where $(u\cdot e_3)\big|_{\interfacep}$ denotes the trace of $u\cdot e_3$ on $\interface$ from above.

\item The level $z=h$ below which we consider the vertical flow to be instantaneous is set such that 
 \begin{equation}\label{coupled_interface}
 \hbot(x) \le  \ha(t,x)\le\max\Bigl\{\min\Bigl\{ \tilde H(t,x)-\frac{\Pbub}{\rho\,g},\hmax(x) \Bigr\},\hbot(x)\Bigr\}, \quad (t,x)\in[0,T[\x\Omega_x.
 \end{equation}

\item The velocity $v$ is defined in $\Omega$ by
\begin{equation}\label{coupled_velocity}
\begin{cases}
 \disp v = u + w &\text{for \ }t\in]0,T[\ ,\quad (x,z)\in\Omega\\
  \disp u = -k_r(P)\,\Big(\frac{1}{\rho\,g}\deriv{P}{z}+1\Big)\,K_0\,e_3 &\text{for \ }t\in]0,T[\ ,\quad (x,z)\in\Omega\\
  \disp w = -k_r\big(  \rho\,g( \tilde H-z )   \big)\,M_0\nabla_x\tilde H& \text{for \ }t\in]0,T[\ ,\quad (x,z)\in\Omega
\end{cases}
\end{equation}

\end{itemize}

 The  coupled model \refmain depends on the definition of the function $h$. Although all intermediate choices respecting \eqref{coupled_interface} are allowed,  we will focus in the next on the two extremal choices
 \begin{align}
  &h(t,x)= \hbot(x),\label{h_hbot}\\
  &h(t,x)= \max\Bigl\{\min\Big\{ \tilde H(t,x)-\frac{\Pbub}{\rho\,g},\hmax(x) \Big\},\hbot(x)\Bigr\}:=h_s(t,x),
 \label{h_H}
 \end{align}
 and on the intermediate one
\begin{equation}\label{h_inter}
h(t,x)= \max\Bigl\{\min\Big\{ \tilde H(t,x)-\frac{\Pbub+R}{\rho\,g},\hmax(x) \Big\},\hbot(x)\Bigr\},
\end{equation}
where $R$ is some positive function possibly depending on $\tilde H$.

\bigskip

 The class of models \refmain is an alternative to the 3D-Richards problem for describing the flow in a shallow aquifer in a large range of time scales. 
This model is designed to fulfill  the two following properties:
\begin{itemize}
 \item to be simpler to handle numerically than the 3D-Richards model
 \item to behave like the 3D-Richards model for any time scale when the  {\it ratio} $\eps$ of the deepness over the horizontal length of the aquifer is small\footnote{however the numerical simulations below show good results even for a {\it ratio} of order 0.1, which is not exceeded by the large majority of the unconfined aquifers.}. For example the behaviors presented in Remark \ref{rem_behavior} are respected. 
\end{itemize}
The first property holds for  \refmain since the 3D original Richards problem is replaced by the coupling of a 2D-problem with a lot of independent 1D-problems which can be solved in parallel. Significant time savings are expected in the computations.
The second property is justified in Section \ref{formal_asymptotic}. The idea is to study the limit $\eps\to0$ of the solution of the 3D-Richards equations and to derive formally the associated effective problem. The same asymptotic analysis is performed for the coupled models \refmain and shows that the corresponding effective problems are exactly the same for every considered time scale and for every choice of $h$ satisfying \eqref{coupled_interface}.

\begin{remark} \label{rem_coupling_time_scales}
It is natural to think that it is possibly not so useful to couple two phenomena which does not hold at the same time scale, since by essence they can not interact with each other.
But the notion of time \textit{scale} is senseless for a  fixed physical situation and we just employ this term to enlighten the interpretations.  
The notion of scale has a precise sense when a sequence of problems is considered, for example parametrised by a small parameter $\eps$ tending to zero with the reference time of study depending on $\eps$. 
This is what we do in Section  \ref{formal_asymptotic} where $\eps$ is the {\it ratio} deepness/length of the aquifer. 
This limit process shows that the two kinds of flow appear at different time scales and then do not interact with each other.
Nevertheless, the coupled problem \refmain \textit{is not} an effective problem and holds without time scale separation assumption. The \textit{depth / width} ratio of the aquifer is then a fixed positive number given by the geometry of the aquifer. 
In particular, "short" and "long" time scales flows can interact without either being negligible or instantaneous.

 \end{remark}

The remainder of this subsection is devoted to comments on the new models \refmain. 
Before splitting those comments according to the choice of the function $h$, we  prove that the model is always mass conservative.

\paragraph{Mass conservation.}
Let $\Mtot(t)$ the total mass of the water contained in domain $\Omega$ at time $t$. 
 We denote by $\Mr$ (resp. $\Md$) the mass of the water filling the domain  $\Omplus$ (resp. $\Ommoins$).
We have 
\begin{gather}
\Mr(t)=\rho\int_{\Omega_x}\int_{\ha(t,x)}^{\hsoil}\phi\,s(P)\,dz\,dx, 
\qquad 
\Md(t)=\rho\int_{\Omega_x}\int_{\hbot(x)}^{\ha(t,x)}\phi\,dz\,dx ,
\label{mass_DR}
\\
\Mtot(t)=\Mr(t)+\Md(t) .
\label{mass_conservation}
\end{gather}
\begin{proposition}\label{prop_mass}
 The total mass satisfies for all  $t\in(0,T)$:
 $$
 \deriv{}{t}\Mtot=-\rho\int_{\Omega_x} (u\cdot e_3)|_{\Gammat}  \, dx .
 $$ 
\end{proposition}
\begin{pf}
By using relation \eqref{mass_DR} and  \eqref{mass_conservation} it comes
\begin{equation}\label{proof_mass_1}
\deriv{}{t}\Mtot=\rho\int_{\Omega_x}\int_{\hbot(x)}^{\ha(t,x)}\phi\deriv{s(P)}{t}\,dz\,dx+\rho\int_{\Omega_x}\int_{\ha(t,x)}^{\hsoil(x)}\phi\deriv{s(P)}{t}\,dz\,dx = \rho\int_{\Omega_x}\int_{\ha(t,x)}^{\hsoil(x)}\phi\deriv{s(P)}{t}\,dz\,dx,
\end{equation}
where the first equality is due to $s(P)=1$ in $]\hbot(x),\ha(t,x)]$ (indeed $P\ge \Pbub$ by \eqref{coupled_watertable} and \eqref{coupled_interface}).
Thanks to the first equation in \eqref{coupled_transition} we deduce
\begin{equation}\label{proof_mass_2}
\int_{\Omega_x}\int_{\ha(t,x)}^{\hsoil(x)}\phi\deriv{s(P)}{t}\,dz\,dx = \int_{\Omega_x} (u\cdot e_3)|_{\interfacep} \,dx-\int_{\Omega_x} (u\cdot e_3)|_{\Gammat}  \,dx.
\end{equation}
Finally by \eqref{coupled_hydraulic} and after an integration by parts
\begin{equation}\label{proof_mass_3}
\int_{\Omega_x} (u\cdot e_3)|_{\interfacep} \,dx=\int_{\partial\Omega_x} \tilde K(\tilde H)\nabla_x\tilde H  \cdot n =0 .
\end{equation}
The result is obtained by plugging \eqref{proof_mass_2} and \eqref{proof_mass_3} in \eqref{proof_mass_1}.
\qed
\end{pf}


\subsection{Comments on the model in the case (\ref{h_hbot})}\label{subsec_comments_bot}


In this case, we have $h=\hbot$, then $\Omplus=\Omega$, $\Ommoins=\emptyset$ and $\interface=\Gammab$ (see \eqref{omplusmoins}). The coupled  model \refmain reduces in: finding the pressure $P$, the velocity $v$ and the auxiliary unknowns $u$, $w$ and $\tilde H$  such that: 
\begin{equation}\label{coupled_velocity_h_hbot}
\begin{cases}
 \disp v = u + w &\text{for \ }t\in]0,T[\ ,\quad (x,z)\in\Omega\\
  \disp u = -k_r(P)\,\Big(\frac{1}{\rho\,g}\deriv{P}{z}+1\Big)\,K_0\,e_3 &\text{for \ }t\in]0,T[\ ,\quad (x,z)\in\Omega\\
  \disp w = -k_r\big(  \rho\,g( \tilde H-z )   \big)\,M_0\nabla_x\tilde H& \text{for \ }t\in]0,T[\ ,\quad (x,z)\in\Omega
\end{cases}
\end{equation}
\begin{equation}\label{coupled_transition_h_hbot}
\begin{cases}
\disp \phi \deriv{s(P)}{t}+\deriv{}{z}\big(u\cdot e_3\big)=0 & \text{for \ }t\in]0,T[\ ,\quad (x,z)\in\Omega\\
 \alpha\, P + \beta\, u\cdot e_3=F  &\text{for \ }(t,x,z)\in]0,T[\x\Gammat\\
  P=\rho\,g\big(\tilde H-\hbot \big)&\text{for \ }(t,x,z)\in]0,T[\x\Gammab\\
  P(0,x,z)=\Pinit(x,z)&\text{for \ } (x,z)\in\Omega\\
\end{cases}
\end{equation}
\begin{equation}\label{coupled_hydraulic_h_hbot}
\begin{cases}
 -\dv_x\big(\tilde K(\tilde H)\,\nabla_x\tilde H \big)=-(u\cdot e_3)|_{\Gammab}&\text{for \ }(t,x)\in]0,T[\x\Omega_x
 \\
 \disp   \tilde K(\tilde H)\,\nabla_x \tilde H\cdot n=0&\text{for \ }(t,x)\in]0,T[\x\partial\Omega_x\\
  \tilde H(0,x)=\Hinit(x)&\text{for \ } x\in\Omega_x
\end{cases}
\end{equation}
This setting corresponds to the simplest form of the model \refmain since \eqref{coupled_hydraulic_h_hbot} is a classical boundary value problem. 
Nevertheless the simulations below illustrate that it is not the better form of approximation for the 3D-Richards equation.

\paragraph{Velocity of the flow.}
The velocity $v$ of the flow turns out to be the superposition of the two velocities $u$ and $w$ which respectively describe the fast and slow components of the flow.
Actually $u$ (resp. $w$) is the dominant component of the flow  in the short time scale (resp. large time scale).

\paragraph{Fast component of the  flow: globally vertical.}
The unknown $u$ represents the velocity associated with the pressure $P$ by the one dimensional Darcy's law given in the second equation of \eqref{coupled_velocity_h_hbot}. This one is deduced from the 3D law (see the second equation of \eqref{richards3d}) by neglecting the horizontal components of the gradient of the pressure $P$.  
By construction the field $u$ is vertical if the conductivity tensor $K_0$ introduced in \eqref{K0} is such that $K_{xz}=0$ but it may admit a non-zero horizontal component in the anisotropic case. 

Furthermore the mass conservation equation \eqref{coupled_transition_h_hbot}  holds.  The pressure $P$ then  satisfies the following vertical Richards equation  where the 
 horizontal variable $x\in \Omega_x$ appears only as a parameter:
\begin{equation}\label{richards_vertical}
 \phi \deriv{s(P)}{t}-\deriv{}{z}\Bigl(k_r(P)\,K_{zz}\,\Big(\frac{1}{\rho\,g}\deriv{P}{z}+1\Big) \Bigr)=0 \qquad\text{in }]0,T[\x\Omega.
\end{equation}
The original 3D-Richards problem reduces to the latter equation when the horizontal diffusion terms are neglected. 
In the short-time scale indeed, those  turn to be non-dominant in  shallow aquifers as  announced in Remark \ref{rem_behavior} and shown in Section \ref{formal_asymptotic}.

The boundary condition on $\Gammat$ remains the same than  in the 3D-Richards problem. 
But on the  bottom  $\Gammab$, the structure of the boundary condition changes and becomes of Dirichlet type, namely $\Pplus\big(t,x,\hbot(t,x)\big)=\rho\,g\big(\tilde H(t,x)-\hbot(t,x) \big)$.  
In fact, even if this Dirichlet condition holds, we do not allow the water  flowing out the aquifer through the bottom boundary. Indeed the possibly non-zero flux $(u\cdot e_3)|_{\Gammab}$ appears  as a source term in the first equation of \eqref{coupled_hydraulic_h_hbot}, so that, as proved in Proposition \ref{mass_conservation}, the coupled  model is globally mass-conservative. 
 The particular value $\Pplus=\rho\,g(\tilde H-\hbot)$ for the bottom Dirichlet condition, has been chosen so that  the fast and slow flows are correctly coupled. 
 This point is further explained in the next paragraph.

\paragraph{Slow component of the flow: globally horizontal.}
On the one hand, introduce the auxiliary pressure $Q$,
$$ Q:=\rho\,g(\tilde H-z),$$
for which $\tilde H$ plays the role of the hydraulic head. 
Since $\tilde H$ does not depend on $z$, we have $(\rho\,g)^{-1}\partial_zQ+1=0$. The first consequence is that the unknown $w$ satisfies  (see \eqref{coupled_velocity_h_hbot})
\begin{equation*}
w=-k_r(Q)\,M_0\nabla_x\tilde H.
\end{equation*}
We recover here the velocity associated to $Q$ by the classical Darcy's law for the conductivity $k_r(Q)\,M_0$. The second consequence is that $Q$ is ruled by
\begin{equation*}
\deriv{}{z}\Big(k_r(Q)\,K_{zz}\,\Big(\frac{1}{\rho\,g}\deriv{Q}{z}+1\Big) \Big)=0 \qquad\text{in }]0,T[\x\Omega ,
\end{equation*}
that is the stationary version of equation \eqref{richards_vertical}. 

On the other hand, we expect $P$ to solve the same stationary problem when the duration of the experiment and when the boundary conditions allow the 1D-Richards problem \eqref{richards_vertical} to reach its stationary state. 
Notice  that such a vertical affine profile is also expected in the 3D-Richards model in any non-short time scale (see Remark \ref{rem_behavior} and Section \ref{formal_asymptotic}).
When this situation occurs, the hydraulic head $H:=P/\rho g+z$ is constant with respect to $z$. The Dirichlet boundary condition on $\hbot$ in \eqref{coupled_transition_h_hbot} then implies that 
\begin{equation*}
H(t,x,z)=H\big(t,x,\hbot(x)\big)=\frac{P\big(t,x,\hbot(x)\big)}{\rho g}+\hbot(x)=\tilde H(t,x).
\end{equation*}
Accordingly, in any non-short time scale, we  get $H\simeq\tilde H$ and then $P\simeq Q$ in $\Omega$. This is the reason of the particular choice $P=\rho\,g(\tilde H-\hbot)$ for the Dirichlet boundary condition on $\hbot$ in \eqref{coupled_transition_h_hbot}.
Roughly speaking, the couple $(Q,w)$ characterizes the flow in a long-time experiment in which the vertical flow seems instantaneous with respect to the horizontal one.

Unlike the velocity $u$, the field $w$ is horizontal both in the isotropic and anisotropic cases due to the  definition of the tensor $M_0$. 
The computations leading to the definition of $M_0$  are done in Section \ref{formal_asymptotic}.
Let us give here some qualitative arguments.
For large times, $w$ is the main order term of the flow which turns out  to be horizontal.
The velocity  $w$ is also related to some hydraulic head, say $L$, by the classical Darcy's law $w=-k_r\,K_0\nabla L$ (as in the Richards equation \eqref{richards3d}; see \eqref{anisotropic_velocity}).
But since $w$ is horizontal we have
\begin{equation*}
0=w\cdot e_3=-k_r\,K_0\nabla L\cdot e_3 = -k_r\,K_{zx}\nabla_x L -k_r\,K_{zz}\deriv{L}{z}\quad \text{and then}\quad\deriv{L}{z}=-k_r\,\frac{K_{zx}}{K_{zz}} \nabla_x L
\end{equation*}
if $K_{zz} \ne 0$ as assumed in this paper, otherwise the question is trivial. 
Accordingly, in the expression of $w=-k_r\,K_0\nabla L$, only the term $\nabla_x L$ appears and it follows $w=-k_r\,M_0\nabla_x L$. Notice that the tensor $M_0$ reduces to  $K_{xx}$ in the isotropic case $K_{xz}=K_{zx}=0$.

Moreover $w$ depends on $z$ only through the term $k_r(\rho\,g(\tilde H-z))M_0$ which decreases to $0$ when $z$ increases above $\tilde H -\Pbub/\rho g$. This decrease is fast in general depending on the soil characteristic function $k_r$. Then, roughly speaking, the horizontal component of the flow is maximum in the saturated part  and  almost  vanishing in the unsaturated one far from the capillary fringe.

 The evolution of the ``stationary pressure'' $Q$ is ruled by the first equation of \eqref{coupled_hydraulic_h_hbot}. This is an horizontal mass-conservation equation associated with the average velocity $\tilde w:=-\tilde K(\tilde H)\nabla_x \tilde H=\int_{\hbot}^{\hsoil}w\,dz$.   The right-hand side is the source term computed from the 1D-Richards problem and  which transfers the mass from the vertical description to the horizontal one.
 
  Notice that in this model \eqref{coupled_velocity_h_hbot}-\eqref{coupled_hydraulic_h_hbot}, the Dupuit hypothesis is not considered. We precise this point in the next Subsection.


\subsection{Comments on the model in the cases (\ref{h_H}) and (\ref{h_inter})}\label{subsec_comments_inter}


Now we come back to the model \refmain in which we set the virtual interface $h$ by
\begin{equation}
h(t,x)= \max\Big\{\min\Big\{ \tilde H(t,x)-\frac{\Pbub+R}{\rho\,g},\hmax(x) \Big\},\hbot(x)\Big\},
\label{h_R_choice}
\end{equation}
for a given non-negative function $R$ possibly depending on $\tilde H$. In the numerical simulations at the end of this section, we consider the constant cases $R=0$, corresponding to  \eqref{h_H}, and $R=3$.
Choosing \eqref{h_H} could be guessed as the most intuitive choice since it means in general splitting the domain along the water table, thus separating the flows in the saturated  and in the unsaturated areas. 
But simulations show that it is not necessary  the optimal choice  for the quality of the 3D-Richards approximation.

\paragraph{Velocity of the flow.}
As previously, the velocity $v$ of the flow results from the contribution of a fast component $u$ and of a slow one $w$.
The set $\Ommoins$ is no more  empty in general and an additional brick is introduced in the model for describing the flow in this area. 
We start by giving some properties of the interface $\interface$.

\paragraph{Interface discriminating the flow behaviors.}

As seen in \eqref{omplusmoins},  the sets $\Ommoins(t)$ and $\Omplus(t)$ are characterised by $\ha$.  In view of  the constraint   \eqref{coupled_interface}, the condition 
\begin{equation}\label{bound.ha}
 \hbot(x) \le \ha(t,x)\le\hmax(x) 
\end{equation}
holds for all $(t,x)\in\domt\x\Omega_x$. 
Due to \eqref{coupled_watertable} and \eqref{coupled_interface} the pressure at the level $z=\ha(t,x)$ satisfies for all $(t,x) \in\domt\x \Omega_x$:
\begin{equation}\label{phazero}
P\big(t,x,\ha(t,x)\big)\begin{cases}
=\Pbub+R\qquad \text{if } \quad \hbot(x)<\ha(t,x)<\hmax(x) ,\\
                        \ge \Pbub+R\qquad \text{if } \quad \ha(t,x)=\hmax(x), \\
			\le \Pbub+R\qquad \text{if } \quad \ha(t,x)=\hbot(x).
                       \end{cases}
\end{equation}
In particular, thanks to \eqref{sat} and since $R\ge0$ we get
\begin{equation}\label{saturated}
s\big(P(t,x,z)\big)=1\qquad \text{if } \quad \hbot(x)< z\le\ha(t,x) ,
\end{equation}
  which means that the set $\Ommoins(t)$ contains a saturated part of the aquifer for any choice of $R\ge0$.  
 More precisely, the soil is fully saturated in $\Ommoins(t)$ for every $t\in\domt$  if  $R>0$, and if $R=0$, that is for \eqref{h_H},  $\Ommoins$ can be interpreted as the water table (see Remark \ref{rem_hsat}   below for precisions).

 By construction  $\ha(t,x)\le\hmax$ so that the interval $]\ha(t,x),\hsoil(x)[$  remains non-empty
for all $(t,x)\in\domt\x\Omega_x$. Then we do not have to explicit a direct coupling of the flow in $\Ommoins$ with the one in the overland. The coupling between $\Ommoins$ and $\Omplus$ is  sufficient.

\paragraph{Fast component of the  flow: globally vertical, a part being instantaneous.}
We start by remarking that, as in the previous case,  the velocity $u$ is related to $P$ by the  vertical Darcy's law \eqref{coupled_velocity}. Moreover the same 1D-Richards equation \eqref{coupled_transition} holds, but now, only in the upper part of the aquifer. In particular, in the short-time scale, the  dominant vertical flow in $\Omplus(t)$ remains well described. 

The main difference between cases $h=\hbot$ and $h\neq \hbot$ is related to the vertical flow in the saturated area $\Ommoins(t)$.   
Indeed,  the pressure profile \eqref{coupled_watertable}  now holds in $\Ommoins$ and  in particular $u$ is zero in $\Ommoins$. As said before, this affine profile is expected in the 
non-short time scale when the vertical flow appears instantaneous. Hence,  the model \refmain describes precisely the vertical flow in $\Omplus$ and assumes that this flow is instantaneous in
$\Ommoins$.  Such an assumption is classical in  models of  saturated  shallow aquifers and is known as the Dupuit hypothesis. Then, the model \refmain in the cases \eqref{h_inter} can be seen as the coupling of a Dupuit horizontal flow in a saturated part at the bottom of the aquifer with many vertical 1D-Richards flows for a precise description of the leaking fluxes from the overland to the water table. 

 Notice that, even if $h\ne\hbot$, the model \refmain does  approximate the 3D-Richards problem at every time scale when the ration $\eps=$\textit{deepness\,/\,horizontal length} tends to zero. Indeed, Proposition \ref{prop_richards} below holds for any choice of function $h$ such that \eqref{coupled_interface} is satisfied. This is explained  by the following points in short times:
\begin{itemize}
 \item From the 3D-Richards problem, we expect a vertical description given by the 1D-Richards in the whole $\Omega$, with a vanishing flux at the bottom of the domain (see \eqref{effective_short}).
 \item From our model, we get 1D-Richards only in $\Omplus$ with a zero flux in $\Ommoins$ (see proof of the short-time scale near equation \eqref{short_rescaled_coupled}) \emph{and} the continuity of the pressure. 
\end{itemize}
In fact,  these problems are exactly the same.

The field $u$ is non-singular thanks to the continuity condition satisfied by $P$ on $\interface$ (see \eqref{coupled_transition} and \eqref{coupled_watertable}). 
As for $h=\hbot$, the particular value of the Dirichlet condition on $\interface$ has been chosen for a proper coupling of the fast and slow components of  the flow. This is further developped in the next paragraph.
However if $u\cdot e_3$ has a trace on the boundary $\interface$ of $\Omplus$, this one is non-zero in general whereas $u\cdot e_3=0$ in $\Ommoins$. This is a notable difference with the case $h=\hbot$.

\paragraph{Slow component of the  flow.}
Again, we introduce  the auxiliary pressure $Q=\rho\,g(\tilde H-z)$  and we remark that now   $P=Q$ in $\Ommoins(t)$ (even for short times). The fact that  $P\simeq Q$ in the whole $\Omega$ for any non-short times  comes, as in the case $h=\hbot$, from the Dirichlet condition $P=\rho\,g(\tilde H-z)$ which holds on $\interface$.

The evolution of $(Q,w)$ is characterized by the evolution of $\tilde H$ given in \eqref{coupled_hydraulic}. In this case where $\Ommoins(t)$ in non-empty in general, we can explicit a little more the dynamic of  $\tilde H$. This is detailed in the next paragraph.

\paragraph{Evolution of the hydraulic head.}
Rewrite the  problem \eqref{coupled_hydraulic} using the first equation of \eqref{coupled_transition} averaged on $[\ha,\hsoil]$:
\begin{equation}
-\dv_x\Big(\tilde K(\tilde H)\,\nabla_x\tilde H \Big) = -u\big|_{\Gammat}\cdot e_3-\int_{h(t,x)}^{\hsoil(x)} \phi\,\deriv{s(\Pplus)}{t} \,dz \qquad\text{in }\domt\x\Omega_x.
\label{alternate_dupuit_-1}
\end{equation}
Since  $s(\Pplus)=1$ for $z\in[\hbot,h]$, we get
\begin{equation}\label{alternate_dupuit_0}
-\dv_x\Big(\tilde K(\tilde H)\,\nabla_x\tilde H \Big) = -u\big|_{\Gammat}\cdot e_3-\deriv{}{t}\int_{\hbot(x)}^{\hsoil(x)} \phi\,s(\Pplus) \,dz \qquad\text{in }\domt\x\Omega_x , 
\end{equation}
 or equivalently by using the Leibniz rule in \eqref{alternate_dupuit_-1} and $s(P)_{\vert z=h}=1$:
\begin{equation}\label{alternate_dupuit}
\phi|_{\interface} \deriv{h}{t}-\dv_x\Big(\tilde K(\tilde H)\,\nabla_x\tilde H \Big) \\
 =-u\big|_{\Gammat}\cdot e_3-\deriv{}{t}\left(\int_{h(t,x)}^{\hsoil(x)} \phi\,s(\Pplus) \,dz \right)\qquad\text{in }\domt\x\Omega_x.
\end{equation}
The hydraulic head $\tilde H$ is characterized by the latter equation completed by the limit conditions in \eqref{coupled_interface}. 
This problem  is a non-linear degenerate diffusion equation. Indeed, the  diffusion tensor $\tilde K(\tilde H)$ vanishes  when $\tilde H$ tends to $-\infty$.
If moreover \eqref{h_H} holds,   in view of \eqref{coupled_interface}, the time derivative can be expressed as 
\begin{equation*}
\deriv{h}{t}=C(\tilde H)\,\deriv{\tilde H}{t}\qquad\text{with}\qquad C(\tilde H)=\begin{cases}
																								    1& \text{if }\tilde H-\Pbub/\rho\, g\in]\hbot,\hmax[\\
																								    0&\text{if not}.
																								    \end{cases}
\end{equation*}

The right-hand side of the first equation in \eqref{coupled_hydraulic} plays the role of a \textit{source term} and represents for each $x\in \Omega$ the evolution of the amount of water which flows in or out the column $]\ha(t,x),\hsoil(x)[$ through its lower boundary $\ha(t,x)$. As we have shown in Proposition \ref{prop_mass} above, this source term ensures the mass conservation in the coupled model \refmain. 
Of course this term also depends (non linearly) on the solution $\tilde H$.
However this dependence is more easy to handle than the one given in the first equation of  \eqref{coupled_hydraulic}. In particular, the expression \eqref{alternate_dupuit} is well adapted to the numerical implementation of the coupled problem \refmain.

Notice that the level $z=h_s$, defined in \eqref{h_H}, represents the interface between the saturated and unsaturated part of the aquifer according to the auxiliary pressure $Q:=\rho\,g(\tilde H(t,x)-z)$. In particular  $Q(t,x,h_s(t,x))=\Pbub$ if  $h_s(t,x)\in(\hbot(x),\hsoil(x))$ (regardless of the choice of $R\ge0$ in \eqref{h_R_choice}).
The conductivity tensor $\tilde K(\tilde H)$ defined in \eqref{averaged_conductivity} can be then decomposed into two parts:
\begin{equation}\label{conductivity_split}
\tilde K(\tilde H)(t,x) =\tilde C_0+ \int_{h_s(t,x)}^{\hsoil(x)}k_r(Q)\,M_0(x,z)\,dz
\end{equation}
where $\tilde C_0$ is the averaged conductivity of the saturated soil, {\it i.e.}
\begin{equation*}
\tilde C_0 =\int_{\hbot(x)}^{h_s(t,x)}M_0(x,z)\,dz.
\end{equation*}
In classical models for the saturated part of an aquifer obtained by  vertical integration under the Dupuit's assumption, the definition of the effective conductivity   (see for example
\cite{Bear1987}) reduces to $\tilde C_0$ instead of  $\tilde K(\tilde H)$, the latter being a little greater. The quantity   
$\tilde C_0$  takes into account the horizontal flow  in the saturated part but it ignores the (little) one  in the unsaturated part, in particular close to the interface $z=h_s$ where the
capillary effects lead to a non-negligible saturation. 
In practice, the smaller $h_s$, the more significant is the difference $\tilde K(\tilde H)-\tilde C_0$. In particular, if   a part of  the bottom of the aquifer is not saturated, that is $h_s=\hbot$,  considering only the vanishing conductivity $\tilde C_0$ whereas $\tilde K(\tilde H)$ remains positive is physically incorrect.

\subsection{Numerical simulations}\label{numeric}

In this section we compare numerically the original 3D-Richards model \eqref{richards3d} and the coupled model \refmain for several choices of $h$ satisfying \eqref{coupled_interface}. 

\begin{figure} 
\begin{center}
     \begin{tikzpicture}[baseline = (a.center)]
 \node[right](a) at (0.6,0.15) {Impermeable rock};  

\fill[color=black!40] (0,0) -- (0.6,0) -- (0.6,0.3) -- (0,0.3) -- cycle; 
 
\end{tikzpicture}
\hfil
     \begin{tikzpicture}[baseline = (a.center)]
 \node[right](a) at (0.6,0.0) {$\hsat$};  
\draw[color=black,dashed,line width=1pt]   (0.0,0) -- (0.48,0) ;  

\end{tikzpicture}
\hfil
     \begin{tikzpicture}[baseline = (a.center)]
 \node[right](a) at (0.6,0)  {$\hsoil$, $\hbot$};  
\draw[color=black!60,line width=1.5]   (0,0) -- (0.5,0) ;  

 \end{tikzpicture}   
\end{center}

 \begin{minipage}{0.495\linewidth}
 \resizebox {\linewidth} {!} {
\begin{tikzpicture}[]
\node[anchor=south west,inner sep=0] (image) at (0,0) {\includegraphics[width=\linewidth,height=0.8\linewidth,trim=4cm 3.cm 3cm 2.2cm,clip]{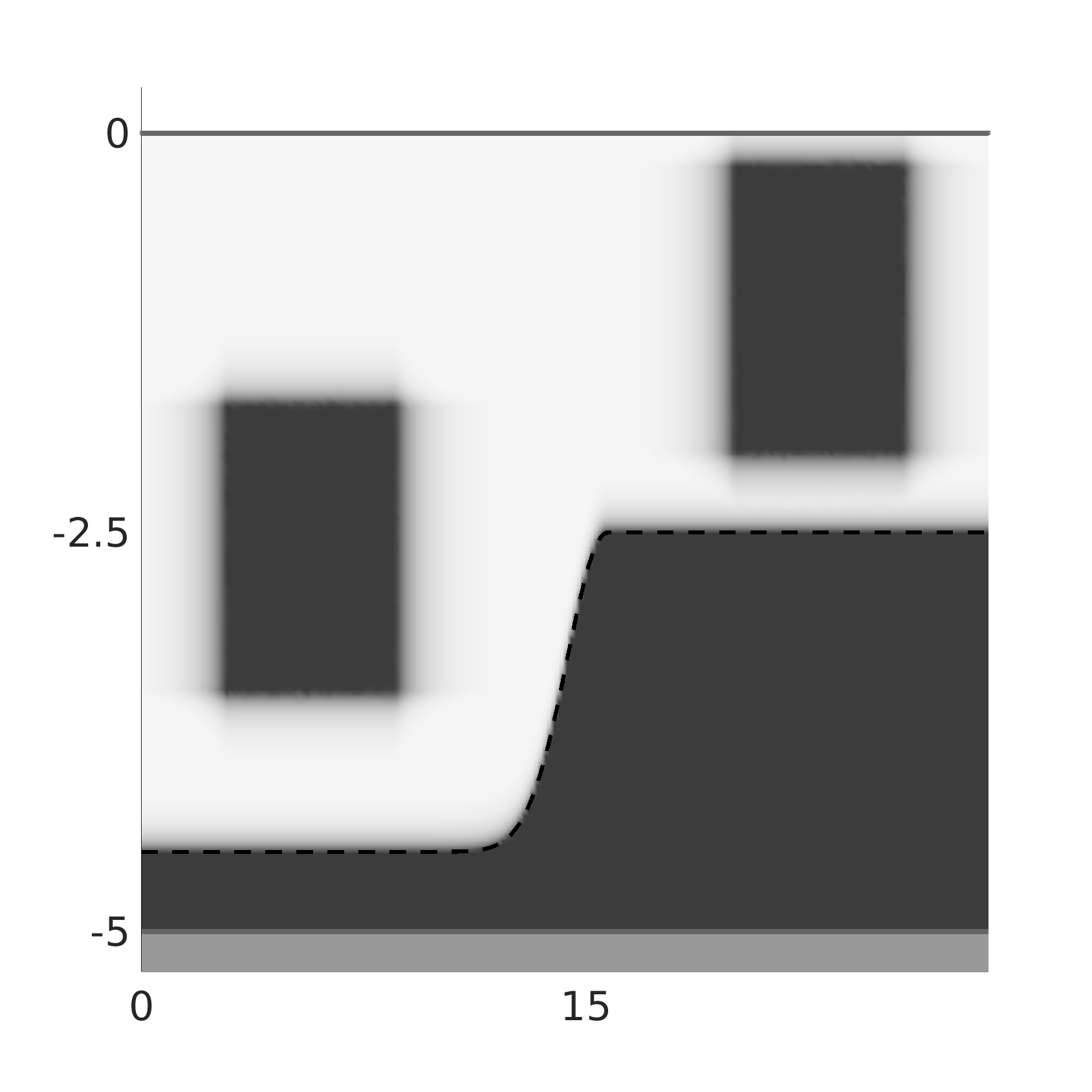}};
    \begin{scope}[x={(image.south east)},y={(image.north west)}]
        \draw (0.2,0.8) node {$t=0$};
        \draw[white] (0,0) node[below] {$0$};
        \draw[white] (0.5,0) node[below] {$x$};                            
        \draw[white] (0.945,0) node[below] {$28$}; 
        \draw (0,0.955) node[left] {$0$};
        \draw (0,0.5) node[left] {$z$};                            
        \draw (0,0.0345) node[left] {$-5$}; 
        \draw[black!20] (0.185,0.47) node {$R_1$};         
        \draw[black!20] (0.755,0.75) node {$R_2$};         

    \end{scope}

\end{tikzpicture}
}
 \end{minipage}
 \hfill\begin{minipage}{0.495\linewidth}
 \resizebox {\linewidth} {!} {
 \begin{tikzpicture}[]
 
 \node[anchor=south west,inner sep=0] (image) at (0,0) {\includegraphics[width=\linewidth,height=0.8\linewidth,trim=4cm 3.cm 3cm 2.2cm,clip]{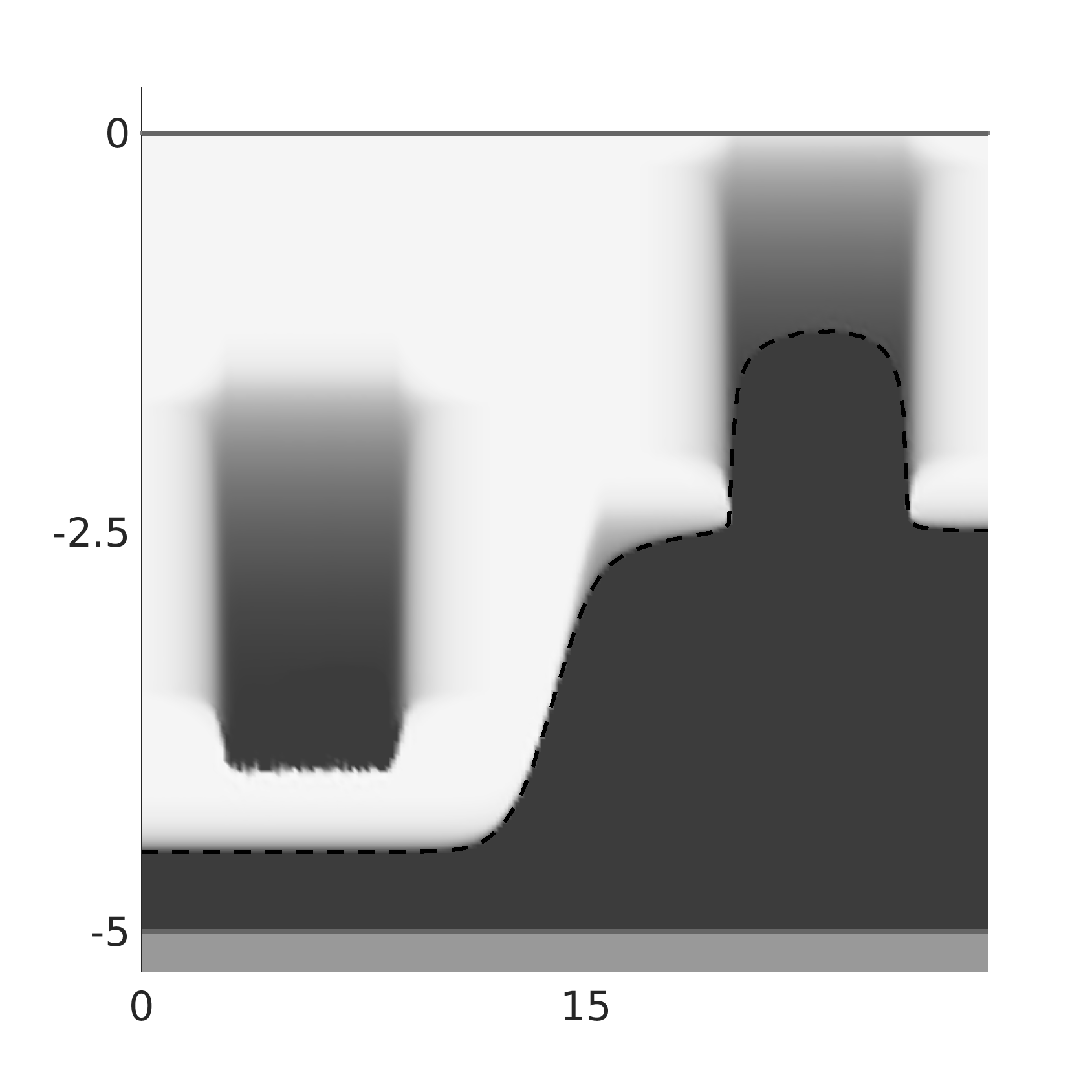}};
    \begin{scope}[x={(image.south east)},y={(image.north west)}]
        \draw (0.2,0.8) node {$t=10$};
        \draw[white] (0,0) node[below] {$0$};
        \draw[white] (0.5,0) node[below] {$x$};                            
        \draw[white] (0.945,0) node[below] {$28$};       
        \draw[white] (0,0.955) node[left] {$0$};
        \draw[white] (0,0.5) node[left] {$z$};                            
        \draw[white] (0,0.0345) node[left] {$-5$}; 

    \end{scope}

\end{tikzpicture} 
 }
 \end{minipage}
 
 \begin{minipage}{0.495\linewidth}
 \resizebox {\linewidth} {!} {
    \begin{tikzpicture}[]
 \node[anchor=south west,inner sep=0] (image) at (0,0) {\includegraphics[width=\linewidth,height=0.8\linewidth,trim=4cm 3.cm 3cm 2.2cm,clip]{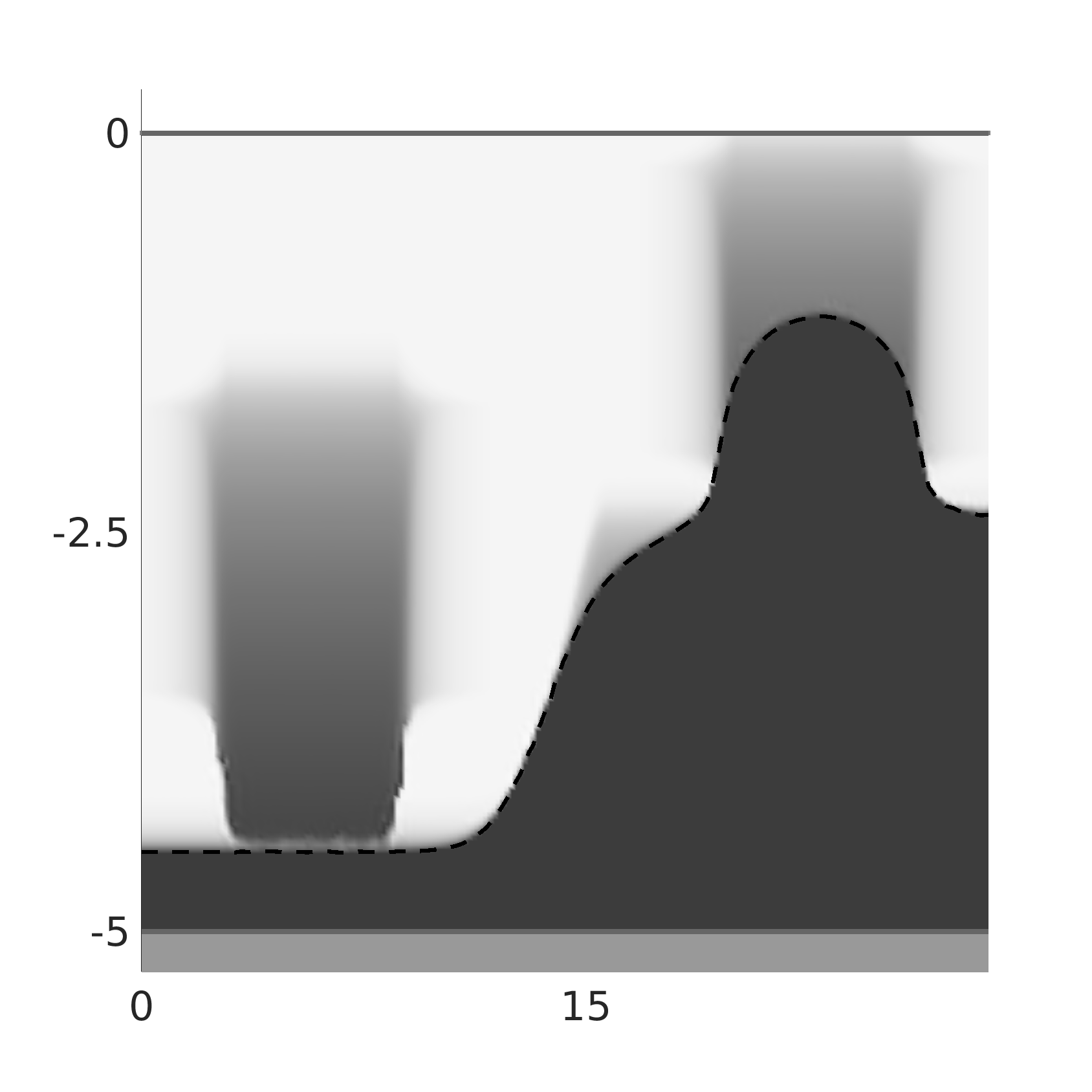}};
    \begin{scope}[x={(image.south east)},y={(image.north west)}]
        \draw (0.2,0.8) node {$t=20$};
        \draw (0,0.955) node[left] {$0$};
        \draw (0,0.5) node[left] {$z$};                            
        \draw (0,0.0345) node[left] {$-5$};                             
        \draw (0,0) node[below] {$0$};
        \draw (0.5,0) node[below] {$x$};                            
        \draw (0.945,0) node[below] {$28$};  

    \end{scope}                  

\end{tikzpicture}
}
 \end{minipage}
 \hfill
 \begin{minipage}{0.495\linewidth}
        \resizebox {\linewidth} {!} {
\begin{tikzpicture}[]
 \node[anchor=south west,inner sep=0] (image) at (0,0) {\includegraphics[width=\linewidth,height=0.8\linewidth,trim=4cm 3.cm 3cm 2.2cm,clip]{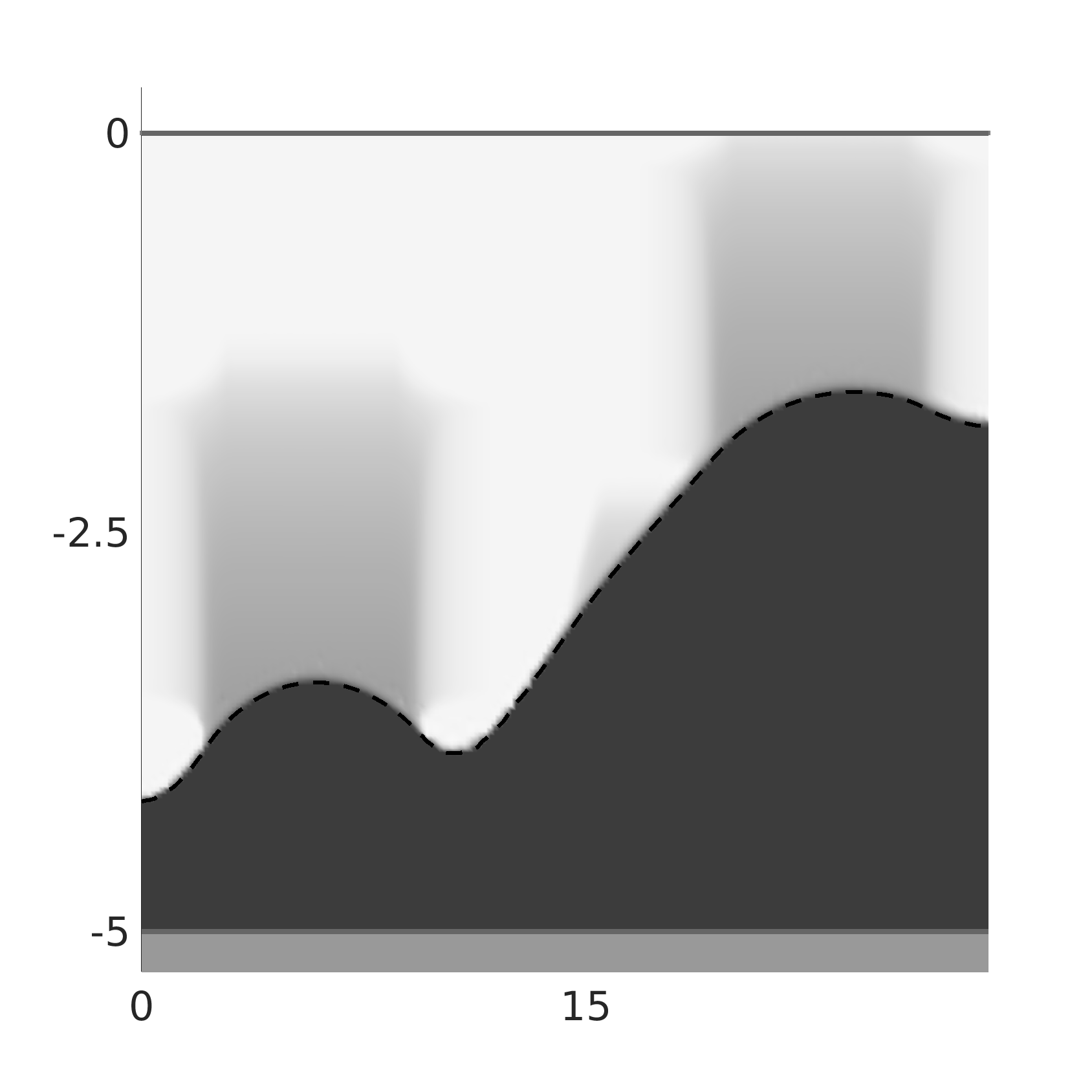}};
    \begin{scope}[x={(image.south east)},y={(image.north west)}]
        \draw (0.2,0.8) node {$t=96$};
        \draw[white] (0,0.955) node[left] {$0$};
        \draw[white] (0,0.5) node[left] {$z$};                            
        \draw[white] (0,0.0345) node[left] {$-5$}; 
        \draw (0,0) node[below] {$0$};
        \draw (0.5,0) node[below] {$x$};                            
        \draw (0.945,0) node[below] {$28$};  

    \end{scope}

\end{tikzpicture}
}
 \end{minipage}

 \begin{center}
 \begin{tikzpicture}[]
\fill[left color=white,right color=black!76]
(0,0) -- (10,0) -- (10,0.5) -- (0,0.5) -- cycle; 
\draw[color=black]
(0,0) -- (10,0) -- (10,0.5) -- (0,0.5) -- cycle;

 \draw (0,0) node[below] {$s=0$};  
 \draw (5,0) node[below] {$s=0.5$};  
 \draw (10,0) node[below] {$s=1$};  
\end{tikzpicture}

\end{center}

\caption{Solution of the classical Richards problem in the reference test case. \label{fig_richards_2d}}
\end{figure}
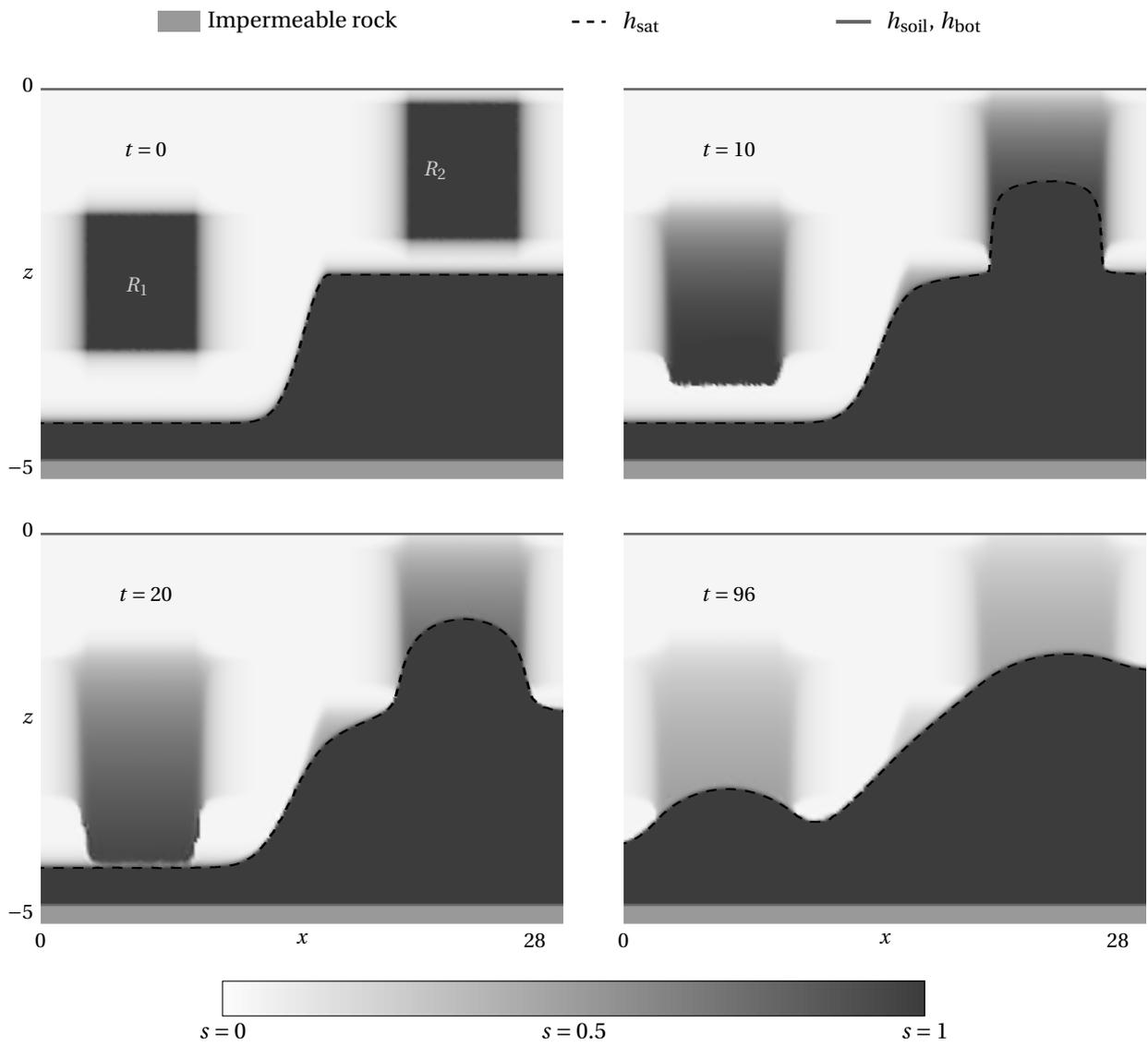

\begin{figure} 
 
 \begin{center}
     \begin{tikzpicture}[baseline = (a.center)]
 \node[right](a) at (0.6,0.15) {Impermeable rock};  

\fill[color=black!40]
(0,0) -- (0.6,0) -- (0.6,0.3) -- (0,0.3) -- cycle; 
 
\end{tikzpicture}
\hfil
     \begin{tikzpicture}[baseline = (a.center)]
 \node[right](a) at (0.6,0.04)  {$\hsoil$, $\hbot$};  
\draw[color=black!60,line width=1.5]   (0,0) -- (0.5,0);

 \end{tikzpicture}

     \begin{tikzpicture}[baseline = (a.center)]
 \node[right](a) at (0.6,0.03) {$\hsatdd$};  

\draw[color=black,line width=1]   (0,0) -- (0.5,0) ;  
\end{tikzpicture}
 \hfil
     \begin{tikzpicture}[baseline = (a.center)]
 \node[right](a) at (0.6,0.04)  {$\hsatmin$};  
\draw[color=black!50,dash dot,line width=1]   (0,0) -- (0.48,0);

 \end{tikzpicture}
    \hfil
    \begin{tikzpicture}[baseline = (a.center)]
 \node[right](a) at (0.6,0)  {$\hsatmax$};  

 \draw[color=black,dotted,line width=1]   (0,0) -- (0.48,0) ;

 \end{tikzpicture}
 \hfil
     \begin{tikzpicture}[baseline = (a.center)]
 \node[right](a) at (0.6,0) {$\hsatint$};  

\draw[color=black!70,dashed,line width=1]   (0,0) -- (0.48,0) ;  
\end{tikzpicture}
 \end{center}

  \begin{minipage}{0.495\linewidth}

  \begin{tikzpicture}[]
\begin{axis}[width=1.05\linewidth,
axis x line = bottom, axis y line = left,xtick={0,28.5},xticklabels={,},ytick={-5,0},yticklabels={$-5$,$0$}]
\addplot[dashed,black!70,line width=1.3pt]  table[x=X,y=coupled_Pbubp2,col sep=comma] {annim_10110.txt};
\addplot[dotted,line width=1.3pt]  table[x=X,y=coupled_classic, col sep=comma] {annim_10110.txt};
\addplot[dash dot,black!50,line width=1.3pt]  table[x=X,y=bot,            col sep=comma] {annim_10110.txt};
\addplot[line width=0.7pt]  table[x=X,y=2d,             col sep=comma] {annim_10110.txt};
\addplot[black!40,fill] coordinates {(0,-5) (28.56,-5)  (28.56,-5.2) (0,-5.2) };
\addplot[dotted,black!60] coordinates {(-0.5,0.4) };
\addplot[dotted,black!60] coordinates {(29.3,-5.4) };

\addplot[black!60, line width=2] coordinates {(0,0) (28.58,0)  };
\addplot[black!60, line width=2] coordinates {(0,-5) (28.58,-5)  };
\end{axis}  
\draw (1,4) node {$t=10$};                            
\draw(-0.3,2.5)  node {$z$} ;

\end{tikzpicture}

 \end{minipage}
 \hfill\begin{minipage}{0.495\linewidth}

  \begin{tikzpicture}[]
\begin{axis}[width=1.05\linewidth,
axis x line = bottom, axis y line = left,xtick={0,28.5},xticklabels={,,},ytick={-5,0},yticklabels={,}]
\addplot[dashed,black!70,line width=1.3pt]  table[x=X,y=coupled_Pbubp2,col sep=comma] {annim_10255.txt};
\addplot[dotted,line width=1.3pt]  table[x=X,y=coupled_classic, col sep=comma] {annim_10255.txt};
\addplot[dash dot,black!50,line width=1.3pt]  table[x=X,y=bot,            col sep=comma] {annim_10255.txt};
\addplot[line width=0.7pt]  table[x=X,y=2d,             col sep=comma] {annim_10255.txt};
\addplot[black!40,fill] coordinates {(0,-5) (28.56,-5)  (28.56,-5.2) (0,-5.2) };
\addplot[dotted,black!60] coordinates {(-0.5,0.4) };
\addplot[dotted,black!60] coordinates {(29.3,-5.4) };

\addplot[black!60, line width=2] coordinates {(0,0) (28.58,0)  };
\addplot[black!60, line width=2] coordinates {(0,-5) (28.58,-5)  };
\end{axis}  
\draw (1,4) node {$t=24$};

\end{tikzpicture}

 \end{minipage}
 
  \begin{minipage}{0.495\linewidth}

  \begin{tikzpicture}[]
\begin{axis}[width=1.05\linewidth,
axis x line = bottom, axis y line = left,xtick={0,28.5},xticklabels={$0$,$28$},ytick={-5,0},yticklabels={$-5$,$0$}]
\addplot[dashed,black!70,line width=1.3pt]  table[x=X,y=coupled_Pbubp2,col sep=comma] {annim_10500.txt};
\addplot[dotted,line width=1.3pt]  table[x=X,y=coupled_classic, col sep=comma] {annim_10500.txt};
\addplot[dash dot,black!50,line width=1.3pt]  table[x=X,y=bot,            col sep=comma] {annim_10500.txt};
\addplot[line width=0.7pt]  table[x=X,y=2d,             col sep=comma] {annim_10500.txt};
\addplot[black!40,fill] coordinates {(0,-5) (28.56,-5)  (28.56,-5.2) (0,-5.2) };
\addplot[dotted,black!60] coordinates {(-0.5,0.4) };
\addplot[dotted,black!60] coordinates {(29.3,-5.4) };

\addplot[black!60, line width=2] coordinates {(0,0) (28.58,0)  };
\addplot[black!60, line width=2] coordinates {(0,-5) (28.58,-5)  };
\end{axis}  
\draw (1,4) node {$t=48$};                            
\draw(3.2,-0.3)  node {$x$} ;
\draw(-0.3,2.5)  node {$z$} ;

\end{tikzpicture}

 \end{minipage}
 \hfill\begin{minipage}{0.495\linewidth}

  \begin{tikzpicture}[]
\begin{axis}[width=1.05\linewidth,
axis x line = bottom, axis y line = left,xtick={0,28.5},xticklabels={$0$,$28$},ytick={-5,0},yticklabels={,}]
\addplot[dashed,black!70,line width=1.3pt]  table[x=X,y=coupled_Pbubp2,col sep=comma] {annim_10998.txt};
\addplot[dotted,line width=1.3pt]  table[x=X,y=coupled_classic, col sep=comma] {annim_10998.txt};
\addplot[dash dot,black!50,line width=1.3pt]  table[x=X,y=bot,            col sep=comma] {annim_10998.txt};
\addplot[line width=0.7pt]  table[x=X,y=2d,             col sep=comma] {annim_10998.txt};
\addplot[black!40,fill] coordinates {(0,-5) (28.56,-5)  (28.56,-5.2) (0,-5.2) };
\addplot[dotted,black!60] coordinates {(-0.5,0.4) };
\addplot[dotted,black!60] coordinates {(29.3,-5.4) };

\addplot[black!60, line width=2] coordinates {(0,0) (28.58,0)  };
\addplot[black!60, line width=2] coordinates {(0,-5) (28.58,-5)  };
\end{axis}  
\draw (1,4) node {$t=96$};                            
\draw(3.2,-0.3)  node {$x$} ;

\end{tikzpicture}

 \end{minipage}
 \caption{Evolution of the iso-pressure $P=\Pbub$ obtained from the classical Richards equation ($\hsatdd$) and from the coupled model for three choices of $h$ given by \eqref{h_hbot}, \eqref{h_H} and \eqref{h_inter} ($\hsat^\kappa$ for $\kappa\in\{a,b,c\}$ respectively). The test case is the one of Figure~\ref{fig_richards_2d}.\label{fig_compare}}

 \end{figure}
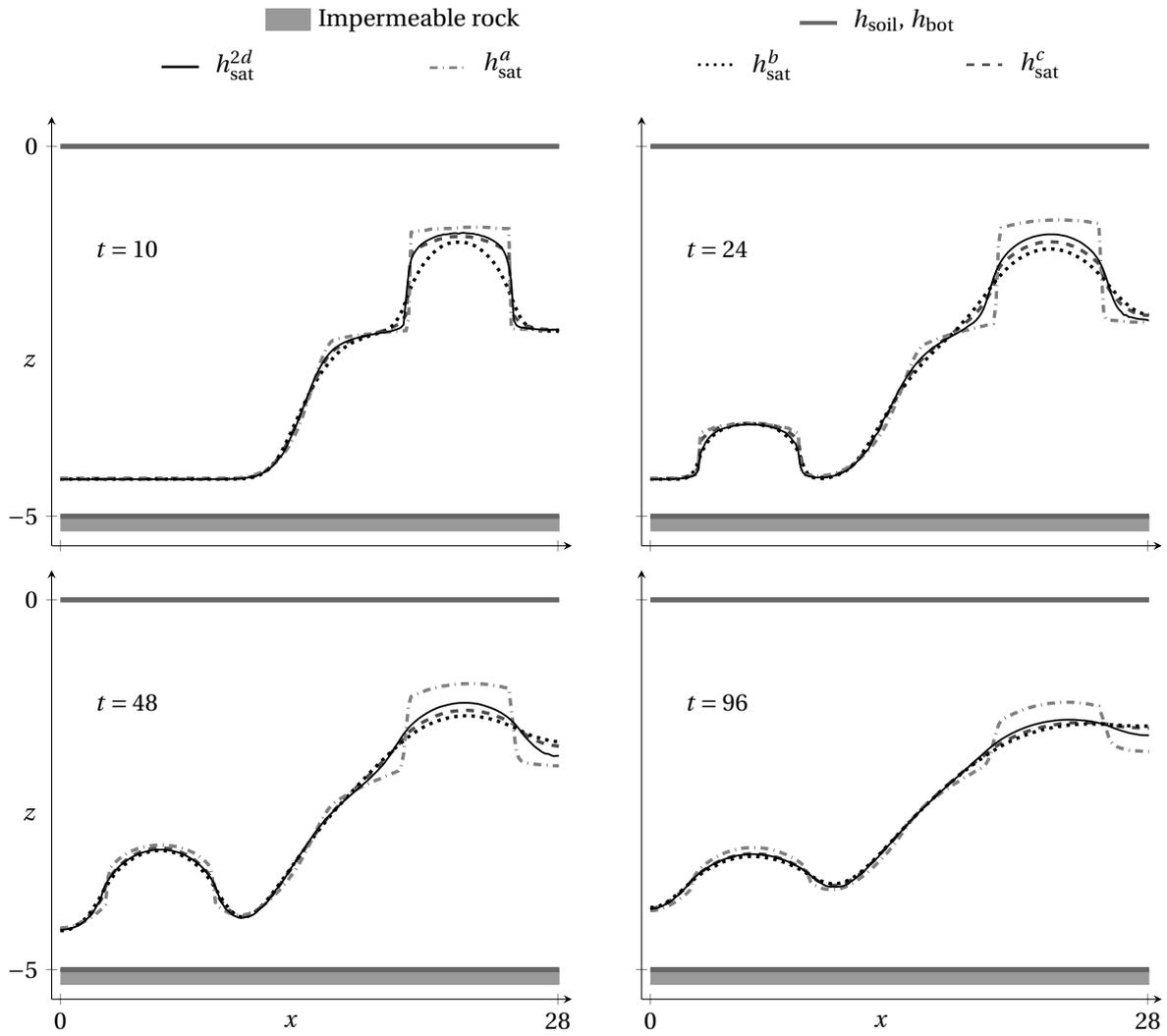

\begin{figure}
     \begin{tikzpicture}[]
\begin{axis}[width=\linewidth,height=0.35\linewidth,legend style={draw=none},legend cell align={left},xtick={0,5.6,10,20,30,40,50,60,70,80},xticklabels={$0$,$5.6$,$10$,$20$,$30$,$40$,$50$,$60$,$70$,$80$},]
   \addplot[mark=x,only marks]    table[x=eps,y=err_bot,col sep=comma] {compar_error.txt}; 
\addlegendentry{\ Minimal case \eqref{h_hbot}}
 \addplot[mark=triangle,only marks]    table[x=eps,y=err_coupled_classic,col sep=comma] {compar_error.txt};
 \addlegendentry{\ Maximal case  \eqref{h_H}}
   \addplot[mark=square,only marks]    table[x=eps,y=err_coupled_Pbubp3,col sep=comma] {compar_error.txt}; 
 \addlegendentry{\ Intermediate case  \eqref{h_inter} for R=3}
\end{axis}
    
\end{tikzpicture}
\caption{Cumulative error in space and time $\|\hsat^{2d}-\hsat^\kappa\|_{L^1([0,T]\times\Omega)}$ {\it versus} the ratio length/deepness$=1/\eps$ of the aquifer ($\kappa\in\{a,b,c\}$). Function $\hsat^{2d}$ is the iso-pressure $P=\Pbub$ in the original 2d-Richards problem and $\hsat^\kappa$ is the one associated with the coupled problems for three different choices of $h$ satisfying \eqref{coupled_interface}.  The test case is the one of figure~\ref{fig_richards_2d}.  \label{fig_error}}

  \end{figure}
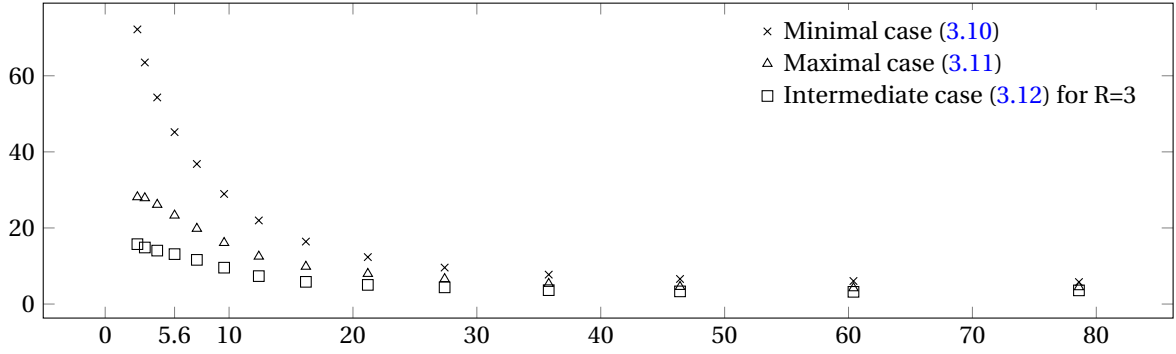

\paragraph{Physical parameters and geometry.}
All the simulations  are done with the following set of data. 
Denoting $I_3$ the $3\x3$ identity matrix  we set:
$$  s(P)=(\Pbub/P)^{\lambda},  \quad  k_r(P)=(\Pbub/P)^{2+3\lambda},   \quad  (\Pbub,\lambda)=(-1.5,3), \quad
  \rho=1 , \quad  \phi=0.1 , \quad K_0=0.1\,I_3 .$$
To lighten the numerical results, we consider the simplified  2D aquifer $\Omega=]-5,0[\x\Omega_x$,    $\Omega_x=]0,L_x[$. 
In the experiments illustrated in Figures \ref{fig_richards_2d} and \ref{fig_compare}, the horizontal length is $L_x=28$. In those of Figure \ref{fig_error},   $L_x \in [21,393]$. 
The parameter $\delta$  in \eqref{hamax} is chosen as small as possible, that is equal to the size of one vertical mesh.
 We assume an impermeable layer at the bottom and the top of the aquifer.

\paragraph{Visualisation.}
For the visualization of the results, we introduce a function $\hsat$ representing in a lot of cases the top level of the saturated region at the bottom of the aquifer ({\it i.e.} the water table).
Let $\hsat=\hsat(t,x)$ and the set $\Ommoinsat(t)$ be defined for a given pressure $P=P(t,x,z)$ by
\begin{equation}\label{hsat}
\hsat(t,x):=\sup I_{t,x},\quad I_{t,x}:=\big\{ z\in[\hbot(x),\hmax(x)]\ |\ P(t,x,z')>\Pbub \ ,\ \forall z'\in[\hbot(x),z[ \big\} ,
\end{equation}
\begin{equation}\label{Omsat}
\Ommoinsat(t):=\bigl\{(x,z)\in\Omega\ |\ z<\hsat(t,x)\bigr\}.
\end{equation}
By construction and if $P$ is continuous we have 
\begin{equation*}
P\big(t,x,\hsat(t,x)\big) \begin{cases}
                             =\Pbub &\text{if }\hbot<\hsat<\hmax\\
                             \ge\Pbub &\text{if }\hsat =\hmax\\
                             \le\Pbub &\text{if }\hsat=\hbot
                           \end{cases}
\end{equation*}
and $P(t,x,z)\ge\Pbub$ for all $z\in]\hbot,\hsat]$. In particular  the soil is fully saturated in $\Ommoinsat(t)$ for every $t\in\domt$.

\begin{remark}\label{rem_hsat}
   Notice that the set $\Ommoinsat$ does not coincide with \textit{the} saturated region of the soil at the bottom of the aquifer. Indeed a saturated region just over $z=\hsat$ is possible for example if $P\ge \Pbub$ also in $\Omega\setminus\Ommoinsat$. The interface $z=\hsat$ then describes
   \begin{itemize}
  \item  either the interface between the saturated part at the bottom of the aquifer and the unsaturated part above in the simplest setting,
  \item or a level between two saturated part when for example a saturated front flow down and reach $\Ommoinsat$, 
  \item or the bottom of the aquifer when $\hsat=\hbot$, that is when there is no saturated part at the bottom,
  \item or the maximum allowed height $\hsat=\hmax$  when, roughly speaking, the water table overflows.
   \end{itemize}
Of course here, since $\hsat(t,x)\le \hsoil-\delta$ by \eqref{hamax},   the set $\Ommoinsat$ cannot reach the soil level $\hsoil$. In this sense  $\Ommoinsat$ does not represent the physical water table which possibly touches the soil level. 
  We only have done this choice for the definition of $\hsat$ to recover the unknown $h$ in the maximal case \eqref{h_H} and thus to facilitate the visualisation. 
\end{remark}

 \paragraph{Numerical scheme.}
For the numerical approximation of the problem \refmain we use mass-conservative fully implicit time schemes associated with finite elements methods in space for both horizontal and vertical directions. The schemes for  \eqref{h_hbot} and \eqref{h_inter} differ: 
\begin{itemize}
 \item In the case \eqref{h_hbot}, we solve directly equation \eqref{coupled_hydraulic} in which the right-hand side $(u \cdot e_3)|_{\interfacep}$ is seen as a Dirichlet to Neumann operator depending on $\tilde H$ and obtained by solving the 1D-vertical Richards equations. This non-linear term is treated with a Newton method.
\item In the case \eqref{h_inter}, the nonlinear coupling between the 1D-vertical Richards equations and the 1D-horizontal diffusion equation is performed by using a Picard's fixed-point method at each time step. This one alternatively solves \eqref{coupled_transition} (for an explicit $\tilde H$ and $h$) and \eqref{alternate_dupuit} (for an explicit right-hand side). 
 \end{itemize}
In any case all the 1D-Richards equations remain  independent at the discrete level and can be solved in parallel.

\paragraph{Reference flowing experiment.}

 At time $t=0$, we consider  a setting where the   function $\hsat$ introduced in \eqref{hsat} corresponds to the height of the water table.
 To show the influence of the deepness of the saturated area, we choose a function $\hsat(0,\cdot)$ which goes smoothly from $-4.5$ on the left part of the aquifer to $-2.5$ on the right one:
  \begin{equation*}
\hsat(0,x)=\begin{cases}
          -4.5 + 2\,e^{-(\frac{15}{L_x})^2(x-0.55\,L_x)^2}&\text{in }[0,0.55L_x],\\
          -2.5&\text{in }]0.55L_x,L_x].
         \end{cases}
\end{equation*}
The initial pressure $P$ is defined by $P(0,x,z)=\rho\,g\,(\hsat(0,x)-z)+\Pbub$ for all $(x,z)$ except near two rectangular regions above $z=\hsat$ where  the pressure goes smoothly to the saturation value $\Pbub$, corresponding to an infiltration process. 
 These rectangles are 
 \begin{equation}\label{rectangles}
R_1=]L_x/10,3L_x/10[\x\left]-3.5,-1.7\right[\and R_2=]7L_x/10,9L_x/10[\x\left]-2,-0.2\right[.
\end{equation}
This initial situation is drawn in the first picture of Figure \ref{fig_richards_2d}. In every picture  the gray scale corresponds to the saturation value, the maximal darkness corresponding to  $s\simeq1$.  
 
The total time of the experiment is $4$ days.  
The solution of the classical Richards problem at time $0$, $10$, $20$ and $96$ hours respectively, is drawn in Figure \ref{fig_richards_2d}. 
The graph of the visualization function $\hsat$ defined in \eqref{hsat} is also plotted. 
Its evolution  will be used for comparing the original Richards model with the coupled model \refmain.  

At  time $t=10$ the water initially in rectangles $R_1$ and $R_2$ started to flow down.
In the right part, some water coming from $R_2$ have reached the saturated water table inducing an increase of its level. In the mean time, we see in the middle of the domain $\Omega_x$  that the water moves to the left and that the function $\hsat$ is  smoother than the initial one.

At  time $t=20$ the water initially in rectangle $R_1$ has continued to flow down and is about to reach the water table. It is important to notice that this flow was essentially along the vertical direction. In particular the water front which is very close to $\hsat$ is approximately horizontal as in the initial situation.

After some time almost all the water initially located  in the rectangle supplies have reached the water table. Then   the interface $\hsat$ becomes flat and is associated with  a pressure admitting the stationary profile $P(t,x,z)=\Pbub+\rho\,g\big(\hsat(t,x)-z\big)$.

\paragraph{Comparison of the models.}

   In this part we compare the solution of the classical Richards model with the one obtained by using the coupled model \refmain. 
We test three particular  choices for the function $h$ satisfying \eqref{coupled_interface}:  the minimal one \eqref{h_hbot}, the  maximal one \eqref{h_H}  and an intermediate one given by \eqref{h_inter} for  $R=3$. 
 All data remain the same as  in the previous paragraph.
In this paper, we focus on the evolution of the functions $\hsat$ defined by  \eqref{hsat}. As indicated in Remark  \ref{rem_hsat}, this function roughly represents the upper level of the water table. In the following we denote by $\hsatdd$ the level  coming from the reference 2d-Richards model and we denote by $\hsatmin$, $\hsatmax$ and $\hsatint$ the ones coming from the model \refmain with the function $h$ given respectively by \eqref{h_hbot}, \eqref{h_H} and \eqref{h_inter}.
   
The functions $\hsatdd$, $\hsatmin$, $\hsatmax$ and $\hsatint$ are plotted in Figure \ref{fig_compare} 
 at time  $t \in \{10,24,48,96\}$ (in hours). 
 We of course do not plot the initial situation which is the same for each model and is the one of the reference test case described in the previous paragraph. The curve $\hsatdd$ is the reference one and is plotted with a black solid line in Figure~\ref{fig_compare}.

Bear in mind that the function $h$ characterizes the level below which the vertical flow is assumed to be instantaneous (instead of being described by the 1D-Richards equation). In every case, the horizontal flow is ruled by equation \eqref{coupled_hydraulic}. 
  \begin{itemize}
   \item In the  case \eqref{h_hbot},  $h=\hbot$. The vertical flow is described by the 1D-Richards model in the whole domain, even in the saturated part   below the level $z=\hsatmin$. 
The horizontal flow in this case seems to be slower than the one given by the Richards model
 (compare the gray dot-dashed line with the black solid one in Figure \ref{fig_compare}).
   \\
   Roughly the idea is that in this case the water have to travel along the whole vertical
   direction before reaching the level $z=h=\hbot$. Then the flux $(u\cdot e_3)|_{\Gammab}$ at the bottom of the aquifer takes a lot of time to increase when the water coming
   from rectangles $R_1$ and $R_2$ reaches  the water table. This flux being the source term in equation \eqref{coupled_hydraulic}, the function $\tilde H$ increases with some delay  
   and the corresponding horizontal flow is slower.

   \item  In the case \eqref{h_H},  $h=\hsatmax$. This case is  opposite of the previous one in the sense that the vertical flow in the whole saturated zone $\Ommoinsat$ is considered to be instantaneous. Then, when the water coming from rectangles $R_1$ and $R_2$ reaches the water table, the flux $(u\cdot e_3)|_{\interface}$ increases very quickly. So does the corresponding hydraulic head $\tilde H$   and the horizontal flow is very and  even too fast  (see the black dotted line compared to the black solid line in Figure \ref{fig_compare}).
   
   \item  In the case \eqref{h_inter} for $R=3$,   $\hbot\le h\le\hsatint$. 
   The corresponding flow should exhibit an intermediate behavior  between the too previous ones.
   Here, the value $R=3$ was chosen so that  $\hsatint$ is very close to the reference one $\hsatdd$ (see the gray dashed line).

  \end{itemize}

  Notice that in every situation, the error between $\hsatdd$ and $\hsat^\kappa$,   $\kappa\in\{a,b,c\}$, is smaller in the left part of the domain than  in the right one. 
  This is due to the fact that  the saturated zone is thiner  in this region. 
  For a very thin saturated region, considering an instantaneous vertical flow or the one given by the vertical 1D-Richards problem gives similar results.   
 Conversely, the thicker the saturated water table is, the more the results issued from the two extremal situations \eqref{h_hbot} and \eqref{h_H} differ from the reference one. 
  Basically,  $\hsat^b$ is expected to move too fast  while  $\hsat^a$  moves too slowly.  
  In this kind of deep situation and if the {\it ratio} between  the deepness and the length of the aquifer is not so small, one of the intermediate choices \eqref{h_inter} is obviously  more appropriate.

   \paragraph{Error made by the coupled model versus the ratio deepness/largeness.}
In the previous simulations, where $\Omega=]0,28[\x]-5,0[$, the {\it ratio} $\eps$=\textit{deepness/length} of the aquifer  is such that $1/\eps=5.6$.  It is important to notice that even in this case of large ratio $\eps$ the error between the original Richards model and the coupled model \refmain in the case \eqref{h_inter} is particularly small (see the dashed plot in Figure \ref{fig_compare}). This supports the fact that the coupled model may be considered for approaching the Richards model also in an aquifer which is not so shallow.
This guess is confirmed by the results plotted in Figure \ref{fig_error}.
The evolution of the error $\|\hsatdd-\hsat^\kappa\|_{L^1(\domt\x\Omega_x)}$ for $\kappa\in \{a,b,c\}$ is drawn in terms of the ratio $1/\eps$. 

As expected all the errors decrease with $\eps$. 
Moreover, the intermediate case \eqref{h_inter} is always the best, mainly in the case of a  ``large'' value of $\eps$. After comes the maximal choice. The worst choice is the maximal one \eqref{h_hbot} but with an error which decreases a lot with $\eps$.

\begin{remark}
 The accuracy of the model depends on the choice of $R$ in \eqref{h_R_choice}, e.g.  for minimizing the error $\|\hsatdd-\hsat^\kappa\|_{L^p(\domt\x\Omega_x)}$. 
This optimization process is postponed to a forthcoming work.
\end{remark}

 \section{Formal asymptotic expansion}\label{formal_asymptotic}


In this section, the 3D-Richards problem \eqref{richards3d} and the coupled model \refmain are compared using asymptotic analysis arguments.
We prove that these models behave the same, whatever the time scale, when the {\it ratio} between the characteristic deepness and the length of the shallow aquifer tends to zero.

\subsection{Dimensionless form of the 3D-Richards and coupled problems}
Introduce a fixed dimensionless reference domain $\ov \Omega$ of type \eqref{omega} and a dimensionless real number $\ov T>0$. 
Fix $\ov\Omega_x$, $\hsoilb$ and $\hbotb$ such that
\begin{equation*}
\ov \Omega=\Big\{(\ov  x,\ov  z)\in\ov \Omega_x\times\R\quad |\quad \ov  z\in\big]\hbotb(\ov  x) ,\hsoilb(\ov x)\big[    \Big\} .
\end{equation*}
To obtain a rescaled version of equations \eqref{richards3d} and \refmain in the domain $]0,\ov T[\x\ov \Omega$, we introduce positive reference numbers $L_x$, $L_z$, $T$. Then, keeping the same notations as in Section \ref{sec_main}, we have:
\begin{itemize}
\item The physical variables are given by
 \begin{equation*}
x =L_x\,\ov x,\quad  z =L_z \,\ov z , \quad  t =  \frac{T}{\ov T}\,\ov t   .
\end{equation*}

\item The corresponding physical domain $\Omega$ is given as in  \eqref{omega} with
\begin{equation*}
\Omega_x = L_x\, \ov\Omega_x,\quad\hsoil(x) =L_z\,\hsoilb(\ov x), \quad  \hbot(x)= L_z\,\hbotb(\ov x).
\end{equation*}
\item  The unknowns are such that
\begin{gather*}
\ov P(\ov t,\ov x, \ov z) = P(t,x,z),\quad \ov v(\ov t,\ov x, \ov z) = v(t,x,z),\quad \ov u(\ov t,\ov x, \ov z) = u(t,x,z),\quad \ov w(\ov t,\ov x, \ov z) = w(t,x,z),
\\
L_z \ov H(\ov t,\ov x) =\tilde H(t,x),\quad L_z \ov h(\ov t,\ov x) =h(t,x) .
\end{gather*}
\item  The reference subdomains are
$$
\Ommoinsb(\ov t)=\big\{(\ov  x,\ov  z)\in\ov \Omega_x\times\R\ |\ \ov  z\in\big]\hbotb(\ov  x) , \hab(\ov t,\ov  x)\big[    \big\},
\qquad
\Omplusb(\ov t)=\big\{(\ov  x,\ov  z)\in\ov \Omega_x\times\R\ |\ \ov  z\in\big]\hab(\ov t,\ov  x) ,\hsoilb(\ov x)\big[    \big\}  .
$$
\item The reference boundaries are $ \Gammabb:=\{( \ov x,\ov z)\in\ov \Omega\ |\ \ov z=\hbotb(\ov x) \}$, $ \Gammatb:=\{(\ov x,\ov z)\in\ov \Omega\ |\ \ov z=\hsoilb(\ov x)\}$ and $ \Gammavb:=\{(\ov x,\ov z)\in\ov \Omega\ |\ \ov  x\in\partial\ov \Omega_x\}$.
\item The reference exterior normals are
\begin{equation*}
\ov n(\ov x,\ov z)=\begin{cases}
  \disp \left(e_3-\frac{L_z}{L_x}\nabla_{\ov x}\hsoilb(\ov x)\right)\left(\frac{L_z^2}{L_x^2}|\nabla_{\ov x}\hsoilb(\ov x)|^2+1\right)^{-1/2}&\text{on }\Gammatb\\
\disp \left(\frac{L_z}{L_x}\nabla_{\ov x}\hbotb(\ov x)-e_3\right)\left(\frac{L_z^2}{L_x^2}|\nabla_{\ov x}\hbotb(\ov x)|^2+1\right)^{-1/2}&\text{on }\Gammabb\\
n(x,z)&\text{on }\Gammavb
\end{cases}
\end{equation*}
where the vector $\ov n$  is horizontal and does not change during the rescaling.
\item  The saturation and relative conductivity satisfy
\begin{equation}\label{ovskr}
s(\ov P) = s(P),\qquad k_r(\ov P) = k_r(P).
\end{equation}
It means that the reference saturation and relative permeability are of order one.
Indeed $P$ and $\ov P$ take the same values, independently of the scale change.
\item  For the conductivities, we set   
\begin{gather}
\ov K_0(\ov x,\ov z)=K_0( x, z),\qquad \ov M_0(\ov x,\ov z)=M_0( x, z),
\label{ovK0}
\\
\ov K(\ov H)(\ov t,\ov x) = L_z\int_{\hbotb(\ov x)}^{\hsoilb(\ov x)} k_r\big(\rho\,g(\ov H(\ov t,\ov x)-\ov z\big)\ov M_0\,d\ov z .
\label{ovKtilde}
\end{gather}
We choose \eqref{ovK0} for the sake of simplicity in the presentation. Indeed, we could also introduce $K$ and $M$ such that $K \ov K_0(\ov x, \ov z)=K_0(x,z)$ and  $M \ov M_0(\ov x, \ov z)=M_0(x,z)$ and then perform the same study assuming that $K/L_x=\mathcal{O}(\eps)$, $M/L_x=\mathcal{O}(\eps)$ and $K/L_z=\mathcal{O}(1)$.
\item  The source term is
\begin{equation*}
\ov F(\ov t,\ov x) = F(t,x) 
\end{equation*}
\end{itemize}

\paragraph{Dimensionless Richards problem.}
Introducing the latter quantities in \eqref{richards3d}, we get the following set of rescaled equations:
\begin{equation}\label{rescaled_mass_richard}
\frac{\ov T}{T}\phi\,\deriv{ s(\ov P)}{\ov t}+\frac1{L_x}\dv_{\ov x}(\ov v) + \frac1{L_z}\deriv{\ov v}{\ov z}=0\quad \text{in }]0,\ov T[\x\ov \Omega ,
\end{equation}
\begin{equation}\label{rescaled_momentum}
 \ov v = -k_r(\ov P)\,\ov K_0\, \Big( \frac1{L_x}\frac{1}{\rho g}\nabla_{\ov x} \ov P + \Big(\frac1{L_z}\frac{1}{\rho g}\deriv{\ov P}{\ov z}+1\big)e_3 \Big)\quad \text{in }]0,\ov T[\x \ov \Omega,
\end{equation}
\begin{equation}\label{rescaled_bottom}
 \ov v\cdot  \Big(\frac{L_z}{L_x}\nabla_{\ov x}\hbotb-e_3\Big) = 0 \quad \text{on }]0,\ov T[\x \Gammabb,
\end{equation}
\begin{equation}\label{rescaled_soil}
\alpha\, \ov P\,\Big(\frac{L_z^2}{L_x^2}\|\nabla_{\ov x}\hsoilb\|^2+1\Big)^{1/2} + \beta\, \ov v\cdot  \Big(e_3-\frac{L_z}{L_x}\nabla_{\ov x}\hsoilb\Big)= \ov F\, \Big(\frac{L_z^2}{L_x^2}\|\nabla_{\ov x}\hsoilb\|^2+1\Big)^{1/2} \quad \text{on }]0,\ov T[\x \Gammatb,
\end{equation}
 \begin{equation}\label{rescaled_vert}
 \ov v\cdot \ov n= 0 \quad \text{on }]0,\ov T[\x \Gammavb .
\end{equation}

Since the aquifer is assumed to be very thin with respect to its horizontal width, the quantity $L_z/L_x$ is very small. 
We choose to consider an aquifer with a fixed height of order $L_z=1$ and a large  horizontal dimension $L_x=1/\eps$ for $\eps\ll1$. We get 
\begin{itemize}
 \item the mass conservation equation which depends on  the time scaling choice $T$:
 \begin{equation}\label{rescaled_mass}
 \disp\frac{\ov T}{T}\phi\,\deriv{ s(\ov P)}{\ov t}+\eps\dv_{\ov x}(\ov v) + \deriv{\ov v\cdot e_3}{\ov z}=0\qquad\text{in }]0,\ov T[\times\ov\Omega
\end{equation}
\item associated with the following Darcy's law and boundary conditions:
\begin{equation}\label{rescaled}
\begin{cases}
  \disp\ov v = -k_r(\ov P)\,\ov K_0\, \left( \frac{\eps}{\rho g}\nabla_{\ov x} \ov P + \Big(\frac{1}{\rho g}\deriv{\ov P}{\ov z}+1\big)e_3 \right)&\text{in } ]0,\ov T[\times\ov\Omega\\
   \disp\alpha\,\ov P\,\Big(\eps^2\,\|\nabla_{\ov x}\hsoilb\|^2+1  \Big)^{1/2} +\beta\,\ov v\cdot  \Big(e_3-\eps\,\nabla_{\ov x}\hsoilb\Big) =\Big(\eps^2\,\|\nabla_{\ov x}\hsoilb\|^2+1  \Big)^{1/2} \, \ov \fluxsoil&\text{on }]0,\ov T[\times\Gammatb\\
 \ov v\cdot \ov n = 0 &\text{on }]0,\ov T[\times\Gammavb \\
 \ov v\cdot  \Big(\eps\nabla_{\ov x}\hbotb-e_3\Big) = 0&\text{on }]0,\ov T[\times\Gammabb
\end{cases}
\end{equation}
\end{itemize}

\paragraph{Dimensionless coupled Dupuit-Richards model.}
By introducing the same parameter $\eps\ll1$, the rescaled coupled problem \refmain reads:
\begin{itemize}
 \item The velocity problem:
\begin{equation}\label{dimensionless_velocity}
\begin{cases}
 \disp \ov v =\ov  u + \ov w &\text{for \ }\ov t\in]0,\ov T[\ ,\quad (\ov x,\ov z)\in\ov \Omega\\
  \disp \ov u = -k_r(\ov P)\,\Big(\frac{1}{\rho\,g}\deriv{\ov P}{\ov z}+1\Big)\, \ov K_0\,e_3 &\text{for \ }\ov t\in]0,\ov T[\ ,\quad (\ov x,\ov z)\in\ov \Omega\\
  \disp \ov w = -\eps\,k_r\big(  \rho\,g( \ov H-\ov z )   \big)\, \ov M_0\nabla_{\ov x}\ov H& \text{for \ }\ov t\in]\ov 0,\ov T[\ ,\quad (\ov x,\ov z)\in\ov \Omega
\end{cases}
\end{equation}
 \item The 1D-Richards equation in the transition zone: 
\begin{equation}\label{dimensionless_transition}
\begin{cases}
\disp \phi\frac{\ov T}{T} \deriv{s(\ov P)}{t}+\deriv{}{\ov z}\big(\ov u\cdot e_3\big)=0 & \text{for \ }\ov t\in]0,\ov T[\ ,\quad (\ov x,\ov z)\in\Omplusb(\ov t)\\
 \alpha\,\ov P+\beta\, \ov u\cdot e_3=\ov F  &\text{for \ }(\ov t,\ov x)\in]0,\ov T[\x\Gammatb\\
  \ov P\big(\ov t,\ov x,\hab(\ov t,\ov x)\big)=\rho\,g\big(\ov H(\ov t,\ov x)-\hab(\ov t,\ov x) \big)&\text{for \ }(\ov t,\ov x)\in]0,\ov T[\x\ov \Omega_x \\
 \ov  P(0,\ov x,\ov z)=\ov \Pinit(\ov x,\ov z)&\text{for \ } (\ov x,\ov z)\in\Omplusb(0)\\
\end{cases}
\end{equation}
\item The  pressure problem in the water table:
\begin{equation}\label{dimensionless_watertable}
       \disp \ov P(\ov t,\ov x,\ov z)=\rho\,g\big( \ov H(\ov t,\ov x)-\ov z \big) \qquad \text{for \ }\ov t\in[0,\ov T[\ ,\quad (\ov x,\ov z)\in\Ommoinsb(\ov t)\\
\end{equation}

\item The hydraulic head problem:
\begin{equation}\label{dimensionless_hydraulic}
\begin{cases}
 \eps^2\,\dv_{\ov x}\Big(\ov K(\ov  H)\,\nabla_{\ov x}\ov H \Big)=\ov u\big|_{\interfacep}\cdot e_3&\text{for \ }(\ov t,\ov x)\in]0,\ov T[\x\ov \Omega_x 
 \\
 \disp  \ov K(\ov  H)  \nabla_{\ov x} \ov H\cdot \ov n=0&\text{for \ }(\ov t,\ov x)\in]0,\ov T[\x\partial\ov \Omega_x \\
 \ov H(0,\ov x)=\ov \Hinit(\ov x)&\text{for \ } \ov x\in\ov \Omega_x 
\end{cases}
\end{equation}
Equivalently, by using  \eqref{alternate_dupuit_0}, the first equation of \eqref{dimensionless_hydraulic}  admits the formulation: for $(\ov t,\ov x)\in]0,\ov T[\x\ov \Omega_x$
\begin{equation}\label{dimensionless_hydraulic_alternate}
 \eps^2\,\dv_{\ov x}\Big(\ov K(\ov  H)\,\nabla_{\ov x}\ov H \Big) \\
 =\ov u\big|_{\Gammatb}\cdot e_3+\frac{\ov T}{T}\deriv{}{\ov t}\left(\int_{\hbotb(\ov x)}^{\ov \hsoil(\ov x)} \phi\,s(\ov \Pplus) \,d\ov z \right)
\end{equation}
\item The definition of the interface separating the two different kind of flows:
\begin{equation}\label{dimensionless_interface}
  \disp   \hbotb(\ov x)\leq  \hab(\ov t,\ov x)\leq\max\left\{\min\Big\{ \ov H(\ov t,\ov x)-\frac{\Pbub}{\rho\,g}, \hmaxb(\ov x) \Big\}, \hbotb(\ov x)\right\}\qquad \text{for \ }(\ov t,\ov x)\in[0,\ov T[\x\ov \Omega_x\\ 
\end{equation}

\end{itemize}


\subsection{Effective problems}


We are interested in the asymptotic behavior of the flow, thus of the models, for both short, intermediate and  large times. 
For the asymptotic analysis, the question is related to the behavior of the dimensionless models above. 
More precisely, we want to describe the effective flow obtained for the short time $T=\ov T$, the intermediate time $T=\eps^{-1}\ov T$ and the long time scales $T=\eps^{-2}\ov T$.

\paragraph{Asymptotic expansion.}
 We introduce the following formal asymptotics for the pressure and the velocity:
\begin{equation}\label{asymptotic_richards}
\disp\ov P_{\eps}^\gamma=\ov  P_0^\gamma+\eps\,\ov  P_1^\gamma+\eps^2\,\ov  P_2^\gamma+\dots\qquad  \ov v_{\eps}^\gamma=\ov  v_0^\gamma+\eps\,\ov  v_1^\gamma+\eps^2\,\ov  v_2^\gamma+\dots
\end{equation}
We emphasize that no arbitrary scaling is imposed, in particular we do not suppose as in \cite{Jazar2014} that the vertical velocity is much smaller than the horizontal one when the ratio $\epsilon$ is very small.
We assume also the existence of formal asymptotics for the  auxiliary variables appearing in \refmain
\begin{equation}\label{asymptotic_coupled}
\begin{cases}
\begin{array}{ll}
\disp\ov u_{\eps}^\gamma=\ov  u_0^\gamma+\eps\,\ov  u_1^\gamma+\eps^2\,\ov  u_2^\gamma+\dots&  \ov w_{\eps}^\gamma=\ov  w_0^\gamma+\eps\,\ov  w_1^\gamma+\eps^2\,\ov  w_2^\gamma+\dots\\
\disp\ov H_{\eps}^\gamma=\ov H_0+\eps\,\ov  H_1^\gamma+\eps^2\,\ov  H_2^\gamma+\dots \qquad&  \disp\ov h_{\eps}^\gamma=\ov h_0^\gamma+\eps\,\ov  h_1^\gamma+\eps^2\,\ov  h_2^\gamma+\dots,
\end{array}
\end{cases}
\end{equation}
and  for the flux at the soil level
\begin{equation}\label{asymptotic_source}
\disp\ov\fluxsoil_{\eps}=\ov  \fluxsoil_0+\eps\,\ov  \fluxsoil_1+\eps^2\,\ov  \fluxsoil_2+\dots.
\end{equation}
Moreover, since $s$ and $k_r$ are $\mathcal{C}^\infty$ by part functions, we write
\begin{equation}\label{asymptotic_sat}
\begin{cases}
\begin{array}{ll}
\disp  s(\ov P_{\eps}^\gamma)= s(\ov  P_0^\gamma)+\eps(\ov  P_1^\gamma+\eps\,\ov  P_2^\gamma+\dots) s'(\ov  P_0^\gamma)+\frac{\eps^2}{2}(\ov  P_1^\gamma+\eps\,\ov  P_2^\gamma+\dots)^2 s''(\ov  P_0^\gamma)+\dots\\
\disp k_r(\ov P_{\eps}^\gamma)= k_r(\ov  P_0^\gamma)+\eps(\ov  P_1^\gamma+\eps\,\ov  P_2^\gamma+\dots) k_r'(\ov  P_0^\gamma)+\frac{\eps^2}{2}(\ov  P_1^\gamma+\eps\,\ov  P_2^\gamma+\dots)^2 k_r''(\ov  P_0^\gamma)+\dots\\
\end{array}
\end{cases}
\end{equation}

\paragraph{Effective problems at the main order.}
Let us introduce the following effective problems:
\begin{itemize}
 \item related to the short time scale ($T=\ov T$),
 \begin{equation}\label{effective_short}
\begin{cases}
\phi\,\disp\deriv{ s(\ov P_0)}{\ov t}+ \deriv{\ov v_0\cdot e_3}{\ov z}=0&\text{in }]0,\ov T[\times\Omega\\
  \disp\ov v_0 = - k_r(\ov P_0) \Big(\frac{1}{\rho g}\deriv{\ov P_0}{\ov z}+1\Big)\ov K_0\,e_3&\text{in }]0,\ov T[\times\Omega \\
\alpha\,\ov P_0+\beta\, \ov v_0\cdot  e_3 = \ov \fluxsoil_0 &\text{on }]0,\ov T[\times\Gammatb\\
 \ov v_0\cdot  e_3 = 0&\text{on }]0,\ov T[\times\Gammabb
\end{cases}
\end{equation}
 \item related to the non-short cases ($T=\eps^{-1}\ov T$ or $T=\eps^{-2}\ov T$ ), 
 \begin{equation}\label{effective_instant}
 \begin{cases}
    \disp\ov P_0(t,x,z)=\rho\,g\big(\ov  H_0(t,x)-\ov  z\big)&\text{in }]0,\ov T[\times\ov \Omega\\
   \disp  \ov  v_0=0&\text{in }]0,\ov T[\times\ov \Omega
   \end{cases}
 \end{equation}
 
  \item related to the non-short cases ($T=\eps^{-1}\ov T$ or $T=\eps^{-2}\ov T$ ) if $\alpha \neq0$
  \begin{equation}\label{effective_alpha_0}
 \ov H_0(\ov t,\ov x)=\frac{\ov F_0(\ov t,\ov x)}{\alpha\,\rho\,g}+\hsoilb(\ov t,\ov x)\qquad\text{in } ]0,\ov T[\times\ov \Omega_x
\end{equation}

 \item related to the intermediate time scale ($T=\eps^{-1}\ov T$)  if $\alpha = 0$ (and then $\beta\neq0$)
\begin{equation}\label{effective_inter}
\rho\,g\Big(\int_{\hbotb}^{\hsoilb} \phi\, s'(\ov P_0)\,dz\Big)\deriv{\ov H_0}{t}  =-\frac{\ov F_1}{\beta}\qquad \text{in } ]0,\ov T[\times\ov \Omega_x
\end{equation}

\item related to the long time scale ($T=\eps^{-2}\ov T$)  if $\alpha = 0$
   \begin{equation}\label{effective_long}
\begin{cases}
 \disp-\dv_x\big(\ov K(\ov H_0)\,\nabla_x\ov H_0\big)=-\frac{\ov  \fluxsoil_2}{\beta}-\deriv{}{\ov t}\Big(\int_{\hbotb}^{\hsoilb}\phi\, s(\ov P_0) \,d\ov z\Big)&\text{in }]0,\ov T[\times\ov \Omega_x\\
\ov K(\ov H_0)\, \nabla_{\ov x}\ov  H_0\cdot \ov n = 0  &\text{on }]0,\ov T[\times\Gammavb
\end{cases}
\end{equation}
and concerning the first order of the velocity
\begin{equation}\label{effective_v1}
\ov  v_1=-\ov k_r(\ov P_0)\,\ov M_0\,\nabla_{\ov x}\ov H_0\qquad\text{in }]0,\ov T[\times\ov \Omega
\end{equation}

\end{itemize}

\begin{proposition}\label{prop_richards}
 Let $(\ov P_{\eps}^\gamma,\ov v_{\eps}^\gamma)$ be the solution of the rescaled 3D-Richards problem \eqref{rescaled_mass}--\eqref{rescaled}
 \emph{or}  of the rescaled coupled model \refrescaled 
  for $T=\eps^{-\gamma}\ov T$ and $\gamma\in\{0,1,2\}$. 
Assume that \eqref{asymptotic_richards}--\eqref{asymptotic_sat} hold true, then  
 \begin{itemize}
  \item $(\ov P_{0}^0,\ov v_{0}^0)$ satisfies \eqref{effective_short}.
  \item $(\ov P_{0}^1,\ov v_{0}^1)$ satisfies  \eqref{effective_instant} and \eqref{effective_alpha_0} if $\alpha\neq0$, or \eqref{effective_instant} and \eqref{effective_inter}  with the compatibility condition $\ov F_0=0$ if $\alpha=0$.
  \item $(\ov P_{0}^2,\ov v_{0}^2)$ satisfies  \eqref{effective_instant} and \eqref{effective_alpha_0} if $\alpha\neq0$, or \eqref{effective_instant} and \eqref{effective_long}  with the compatibility condition $\ov F_0=\ov F_1=0$ if $\alpha=0$.
     Moreover $\ov v_1^2$ satisfies \eqref{effective_v1}  if $\alpha=0$.
\end{itemize}

\end{proposition}

We emphasize  that the intermediate variable $\ov h$ which characterizes the coupled model \eqref{dimensionless_velocity}-\eqref{dimensionless_hydraulic_alternate} does not appear in any of the main order effective problems \eqref{effective_short}-\eqref{effective_long}. 
This agrees with the fact that the whole class of models given by \refmain for \textit{any} $h$ satisfying \eqref{coupled_transition} can approximate the reference Richards model.


\subsection{Proof of Proposition \ref{prop_richards} for the Richards model}


The proof of Proposition \ref{prop_richards} consists in substituting the formal asymptotic expansion \eqref{asymptotic_richards}--\eqref{asymptotic_sat} in the rescaled 3D-Richards problem \eqref{rescaled_mass}--\eqref{rescaled}.
A cascade of equations follows by identifying the powers of $\eps$.
Then we characterize  the main order  terms in the expansion \eqref{asymptotic_richards}. 
In order to reduce ratings in this section, we do not write the exponent $\gamma$ on the variables name.

 \paragraph{General relations.}
Let us start by obtaining the first relations holding in every time scale (i.e. for all $\gamma\in\{0,1,2\}$). 
By plugging the asymptotic expansion \eqref{asymptotic_richards} in the first equation of \eqref{rescaled} we get the following relations holding in $]0,\ov T[\times\Omega$ 
\begin{equation}\label{velocity_gen}
\begin{cases}
 \disp \ov v_0=-k_r(\ov P_0)\,\Big(\frac{1}{\rho\,g}\deriv{\ov P_0}{\ov z}+1\Big)\,\ov K_0\,e_3 ,\\
  \disp \ov v_1=-\frac{k_r(\ov P_0)}{\rho\,g}\,\ov K_0\,\Big(\nabla_{\ov x}\ov P_0+\deriv{\ov P_1}{\ov z }\,e_3\Big)-k_r'(\ov P_0)\,\ov P_1\,\Big(\frac{1}{\rho\,g}\deriv{\ov P_0}{\ov z}+1\Big)\,\ov K_0\,e_3  . 
 \end{cases}
\end{equation}
The same process in the three last equations of \eqref{rescaled} yields the following relations in  $]0,\ov T[$:
\begin{itemize}
 \item on $\Gammatb$
\begin{equation}\label{boundary_gen_soil}
 \begin{cases}
\disp \alpha\,\ov P_0+\beta\,  \ov v_0\cdot  e_3 = \ov\fluxsoil_0 , \qquad
\disp \alpha\,\ov P_1+\beta\, \big( \ov v_1\cdot  e_3-\ov v_0\cdot  \nabla_x\hsoilb\big) = \ov\fluxsoil_1 ,\\
\disp  \alpha\,\left( \ov P_2+\frac12\|\nabla_x\hsoilb\|^2 \,\ov P_0\right)+\beta\, \big(\ov v_2\cdot  e_3- \ov v_1\cdot  \nabla_x\hsoilb\big)= \frac12\|\nabla_x\hsoilb\|^2\,\ov\fluxsoil_0+\ov\fluxsoil_2 ;
 \end{cases}  
 \end{equation}
  \item on $\Gammabb$, for all $k\in \N^* $
\begin{equation}\label{boundary_gen_bottom}
 \disp    \ov v_0\cdot  e_3 =0 ,\qquad 
 \disp \ov v_{k-1}\cdot  \nabla_{\ov x}\hbotb=\ov v_{k}\cdot e_3  ;
\end{equation}
 \item on $\Gammavb$, for all $k\in \N $
\begin{equation}\label{boundary_gen_ver}
 \ov v_k\cdot \ov n = 0 .
\end{equation}

\end{itemize}

\paragraph{Short time case.}
We prove the first claim of Proposition \ref{prop_richards} which is associated with the short characteristic time scale $T=\eps^{-\gamma}\ov T$ for $\gamma=0$. 
The equation \eqref{rescaled_mass} here reads
\begin{equation}\label{rescaled_short}
\phi\,\disp\deriv{s(\ov P)}{\ov t}+\eps\dv_{\ov x}(\ov v) + \deriv{\ov v\cdot e_3}{\ov z}=0 .
\end{equation}
Some computations show that the main order terms in the latter equation give
\begin{equation}
\phi\,\disp\deriv{ s(\ov P_0)}{\ov t}+ \deriv{\ov v_0\cdot e_3}{\ov z}=0\qquad\text{in } ]0,\ov T[ \times \ov\Omega .
\end{equation}
The latter equation completed with  the first equations of \eqref{velocity_gen}, \eqref{boundary_gen_soil}  and \eqref{boundary_gen_bottom} gives exactly the system \eqref{effective_short}. The first claim of Proposition \ref{prop_richards} is proven.

\paragraph{Intermediate time case.}
In this part, we prove the second claim of Proposition \ref{prop_richards} which is associated with the intermediate time scale $T=\eps^{-\gamma}\ov T$ for $\gamma=1$.
Equation of \eqref{rescaled_mass} is now
\begin{equation}\label{rescaled_inter}
\eps\phi\,\disp\deriv{ s(\ov P)}{\ov t}+\eps\dv_{\ov x}(\ov v) + \deriv{\ov v\cdot e_3}{\ov z}=0 .
\end{equation}
We introduce the asymptotic expansion \eqref{asymptotic_richards} in the previous equation and we identify the main order terms. We obtain 
\begin{equation}\label{inter_order_0}
\disp\ \deriv{\ov v_0\cdot e_3}{\ov z}=0\qquad \text{on }]0,\ov T[\times\ov \Omega.
\end{equation}
This constant vertical velocity is actually zero  due to \eqref{boundary_gen_bottom}. 
Moreover, with the first equation of \eqref{velocity_gen} and since $k_r$ and $(\ov K_0)_{33}$ are non-vanishing ($\ov K_0$ is positive definite), we get in $]0,\ov T[\x\ov\Omega$
\begin{equation}\label{velocity_0}
\deriv{\ov P_0}{\ov z}+\rho\,g=0 \qand \ov v_0=0.
\end{equation}
The existence of $\ov  H_0=\ov  H_0(t,x)$ such that
 \begin{equation}\label{P0order0inter}
 \ov P_0(t,x,z)=\rho\,g\big(\ov  H_0(t,x)-\ov  z\big)\qquad\text{in } ]0,\ov T[\x\ov\Omega
\end{equation}
follows. 
Next, since $\ov v_0=0$, the first equation of \eqref{boundary_gen_soil} is 
\begin{equation}\label{P0F0_order1_richards}
\alpha\,\ov P_0=\ov F_0\qquad \text{on }\Gammatb.
\end{equation} 
We now have to differentiate the computations depending on whether  $\alpha=0$  or not.

 If $\alpha\neq0$, then for all $(\ov t,\ov x)\in ]0,\ov T[\x\ov\Omega_x$ we have $\ov P_0(\ov t,\ov x,\hsoilb(\ov t,\ov x))=\ov F_0(\ov t,\ov x)/\alpha$. Accordingly, thanks to \eqref{P0order0inter}, it holds
 \begin{equation*}
\ov H_0(\ov t,\ov x)=\frac{\ov F_0(\ov t,\ov x)}{\alpha\,\rho\,g}+\hsoilb(\ov t,\ov x).
\end{equation*}
 This ends the proof of the second claim of Proposition \ref{prop_richards} in the case $\alpha\neq0$.

 If $\alpha=0$ (then $\beta\neq0$),  equation \eqref{P0F0_order1_richards} only implies that $\ov\fluxsoil_0=0$. In particular, $\ov H_0$ remains as a degree of freedom and we have to exploit the next order terms in the asymptotic expansion for the closure of the effective problem. 
 Identifying the coefficients associated with $\eps^1$  in equation  \eqref{rescaled_inter} we have
 \begin{equation}\label{inter_order_1}
\disp \phi\,\deriv{ s(\ov P_0)}{\ov t}+ \deriv{\ov v_1\cdot e_3}{\ov z}=0\qquad\text{in } ]0,\ov T[ \times \ov\Omega.
\end{equation}
To eliminate $\ov v_1$, we integrate vertically on $]\hbotb,\hsoilb[$ the equation above. After using the fact that $\partial_t (s(\ov P_0))=\rho\,g\,s'(\ov P_0)\,\partial_t \ov H_0$ (consequence of \eqref{P0order0inter}) we have
 \begin{equation}\label{inter_order_11}
\rho\,g\Big(\int_{\hbotb}^{\hsoilb} \phi\, s'(\ov P_0)\,dz\Big)\deriv{\ov H_0}{t} +(\ov v_1|_{\hsoilb}-\ov v_1|_{\hbotb})\cdot e_3 =0.
\end{equation}
Thanks to the second equations of \eqref{boundary_gen_soil} and \eqref{boundary_gen_bottom} in the case where $\alpha=0$ and $\ov v_0=0$, it follows:
\begin{equation*}
\ov v_1\cdot e_3 = \ov \fluxsoil_1/\beta\quad \text{on } \Gammatb \qand \ov v_1\cdot e_3= 0\quad \text{on } \Gammabb.
\end{equation*}
Accordingly, equation \eqref{inter_order_11} becomes
\begin{equation}\label{inter_order_12}
\rho\,g\Big(\int_{\hbotb}^{\hsoilb} \phi\, s'(\ov P_0)\,dz\Big)\deriv{\ov H_0}{t}  =-\frac{\ov F_1}{\beta}.
\end{equation}
 Finally, collecting equations \eqref{P0order0inter} and\eqref{inter_order_12} we get $\ov v_0=0$ and
\begin{equation}\label{inter_final}
\begin{cases}
  \disp\ov P_0(\ov t,\ov x,\ov z)=\rho\,g\big(\ov  H_0(\ov t,\ov x)-\ov  z\big)&\text{in }]0,\ov T[\times\ov \Omega\\
 \disp\rho\,g\Big(\int_{\hbotb}^{\hsoilb} \phi\, s'(\ov P_0)\,dz\Big)\deriv{\ov H_0}{t}  =-\frac{\ov F_1}{\beta} 
 &\text{in }]0,\ov T[\times\ov \Omega_x
\end{cases}
\end{equation}
which correspond to the second claim of Proposition \ref{prop_richards} in the case $\alpha=0$.

 \paragraph{Long time case.}
In this part, we prove the third claim of Proposition \ref{prop_richards} which is associated with the intermediate time scale $T=\eps^{-\gamma}\ov T$ for $\gamma=2$.
Equation  \eqref{rescaled_mass} takes the form
\begin{equation}\label{rescaled_long}
\eps^2\phi\,\disp\deriv{s(\ov P)}{\ov t}+\eps\dv_{\ov x}(\ov v) + \deriv{\ov v\cdot e_3}{\ov z}=0 .
\end{equation}
We substitute the asymptotic expansion\eqref{asymptotic_richards} in the previous equation. The main order part of the equation is
 $\partial_z(\ov v_0\cdot e_3)=0$ which leads, as before, to \eqref{effective_instant} for some function $\ov H_0$ which does not depends on $\ov z$. The same relation \eqref{P0F0_order1_richards} holds and
the characterization of $\ov H_0$ depends on the values of $\alpha$. As before, if $\alpha\neq0$ we have \eqref{effective_alpha_0}.

It remains to deal with the case $\alpha=0$ and to exhibit the equations of system \eqref{effective_long}.
In this case,  the compatibility condition $F_0=0$ holds as before because of  \eqref{P0F0_order1_richards}. 
The characterization of $\ov H_0$ needs to go at the next order in the asymptotic expansion. 
In equation  \eqref{rescaled_long} we get
\begin{equation}\label{long_order_1}
0=\dv_{\ov x}(\ov v_0) + \deriv{\ov v_1\cdot e_3}{\ov z}=\deriv{\ov v_1\cdot e_3}{\ov z}
\end{equation}
where the second equality is due to $\ov v_0=0$. Moreover, the second equations of \eqref{boundary_gen_soil} and \eqref{boundary_gen_bottom} for $k=1$ lead to (since $\alpha=0$)
\begin{equation}\label{long_order_1bis}
 \beta\,\ov v_1\cdot e_3 = \ov \fluxsoil_1\quad\text{on }\Gammatb\qand \ov v_1\cdot e_3= 0\quad\text{on }\Gammabb.
\end{equation}
Then, the  vertical component of the velocity (which is constant by \eqref{long_order_1}) $\ov v_1\cdot e_3$ is zero. Moreover the second compatibility condition $\ov F_1=0$ holds  true thanks to  \eqref{long_order_1bis}.
Using the second equation of \eqref{velocity_gen} and bearing in mind that $(\rho\,g)^{-1}\partial_z \ov P_0+1=0$, we obtain
\begin{equation}\label{anisotropic_velocity}
\ov v_1=-\frac{k_r(\ov P_0)}{\rho\,g}\,\ov K_0\,\Big(\nabla_{\ov x}\ov P_0+\deriv{\ov P_1}{\ov z}\,e_3\Big).
\end{equation}
Since $\ov v_1\cdot e_3=0$, using the same notation for $\ov K_0$ than in \eqref{K0}, we compute  $\partial_z \ov P_1$ by
\begin{equation*}
\deriv{\ov P_1}{\ov z} = -\frac{1}{\ov K_{zz}}\ov K_0\nabla_{\ov x}\ov P_0\cdot e_3.
\end{equation*}
Finally, substitution in the equation above with the relation $\ov P_0=\rho\,g(\ov H_0-z)$ give
\begin{equation}\label{long_order_1final}
\ov v_1=- k_r(\ov P_0)\,\ov M_0\,\nabla_{\ov x}\ov  H_0\quad\text{with}\quad \ov M_0=\begin{pmatrix}
                                                                                        I_2 &-\frac{\ov K_{xz}}{\ov K_{zz}}\\
                                                                                        0&0
                                                                                       \end{pmatrix}\ov K_0=\begin{pmatrix}
                                                                                       \ov S_0&0\\
                                                                                       0&0
                                                                                   \end{pmatrix}     
\end{equation}
where $I_2$ is the $2d$ identity matrix and $\ov S_0=\ov K_{xx}-K_{zz}^{-1}\ov K_{xz}\ov K_{zx}$.

On the other hand, the equation \eqref{boundary_gen_ver} for $k=1$ leads to $\ov v_1\cdot\ov n=0$ on $\Gammavb$. Since $ k_r(\ov P_0)$ does not vanish, we obtain the last equation of \eqref{effective_long}.
After identifying the coefficients associated with $\eps^2$  in equation \eqref{rescaled_long} we get
\begin{equation}\label{long_order_2}
\disp\phi\,\deriv{ s(\ov P_0)}{\ov t}+\dv_{\ov x}(\ov v_1) + \deriv{\ov v_2\cdot e_3}{\ov z}=0.
\end{equation}
Taking into account \eqref{effective_instant}, \eqref{long_order_1final} and the fact that $\alpha=F_0=0$, the third equation of \eqref{boundary_gen_soil} and the second equations of  \eqref{boundary_gen_bottom} for $k=2$ become
\begin{equation}\label{long_order_2boundary}
 \disp \ov v_2\cdot e_3-\ov v_1\cdot  \nabla_{\ov x}\hsoilb = \ov \fluxsoil_2/\beta ,\qquad
 \disp  \ov v_2\cdot e_3- \ov v_1\cdot \nabla_{\ov x}\hbotb = 0 \qquad \text{on }\Gammabb .
\end{equation}
To eliminate $v_2$ in system \eqref{long_order_2}--\eqref{long_order_2boundary}, we integrate \eqref{long_order_2} with respect to $\ov z$ on $[\hbotb,\hsoilb]$. Taking into account the boundary conditions on $\Gammabb$ and $\Gammatb$ we obtain
\begin{equation*}
\deriv{}{\ov t}\int_{\hbotb}^{\hsoilb}\phi\, s(\ov P_0) \,d\ov z+ \int_{\hbotb}^{\hsoilb}\dv_{\ov x} \ov v_1\,d\ov z+\ov v_1|_{\hsoilb}\cdot \nabla_{\ov x}\hsoilb+\frac{\ov \fluxsoil_2}{\beta}-\ov v_1|_{\hbotb}\cdot \nabla_{\ov x}\hbotb=0.
\end{equation*}
We use the Leibniz rule in the second integral and we get
\begin{equation}\label{long_dupuit}
\deriv{}{\ov t}\int_{\hbotb}^{\hsoilb}\phi\, s(\ov P_0) \,d\ov z+ \dv_{\ov x} \Big(\int_{\hbotb}^{\hsoilb}\ov v_1\,d\ov z\Big) = -\frac{\ov F_2}{\beta}.
\end{equation}
Using the averaged conductivity $\ov K$ defined in \eqref{ovKtilde},  we get, with  the first equation of \eqref{long_order_1final}, 
\begin{equation*}
\int_{\hbotb}^{\hsoilb}\ov v_1\,d\ov z = -\int_{\hbotb}^{\hsoilb}k_r(\ov P_0)\, \ov M_0\,\nabla_{\ov x}\ov  H_0=-\ov K(\ov H_0)\,\nabla_{\ov x}\ov  H_0.
\end{equation*}
The above equation associated with  equation \eqref{long_dupuit}  is exactly the system \eqref{effective_long}.
This ends the proof of the last claim of Proposition \eqref{prop_richards}.


\subsection{Proof of Proposition \ref{prop_richards} for the coupled models}


The strategy of the proof  is exactly the same than in the previous subsection.

 \paragraph{General relations.}
Let  $\gamma\in\{0,1,2\}$. 
Using the  expansion \eqref{asymptotic_richards}--\eqref{asymptotic_sat}, we identify powers  of $\eps$ in  all the equations in \eqref{dimensionless_velocity}--\eqref{dimensionless_interface} that does not depend on the time scale $T$.
We obtain from the second  equation of \eqref{dimensionless_velocity} 
\begin{equation}\label{velocity_general1}
\begin{cases}
  \disp \ov u_0 = -k_r(\ov P_0)\,\Big(\frac{1}{\rho\,g}\deriv{\ov P_0}{\ov z}+1\Big)\,\ov K_0\,e_3&\text{in }]0,\ov T[\x \ov \Omega ,\\
  \disp \ov u_1= -\frac{k_r(\ov P_0)}{\rho\,g}\,\deriv{\ov P_1}{\ov z}\,\ov K_0\,e_3-k_r'(\ov P_0)\,\ov P_1\,\left(\frac{1}{\rho\,g}\deriv{\ov P_0}{\ov z}+1\right)\,\ov K_0\,e_3  &\text{in }]0,\ov T[\x \ov \Omega , 
\end{cases}
\end{equation}
from the third equation of \eqref{dimensionless_velocity} 
\begin{equation}\label{velocity_general2}
  \disp \ov w_0=0 ,\qquad
  \disp \ov w_1=-k_r\big(\rho\,g(\ov H_0-\ov z)\big) \ov M_0\,\nabla_{\ov x} \ov H_0 \qquad \text{in }]0,\ov T[\x \ov \Omega ,
\end{equation}
from the first equation of \eqref{dimensionless_velocity} 
\begin{equation}\label{velocity_global}
\begin{cases}
    \disp       \ov v_0= \ov u_0+ \ov w_0= \ov u_0=-k_r(\ov P_0)\,\Big(\frac{1}{\rho\,g}\deriv{\ov P_0}{\ov z}+1\Big)\,\ov K_0\,e_3&\text{in }]0,\ov T[\x \ov \Omega ,\\ 
    \disp        \ov v_1= \ov u_1+ \ov w_1&\text{in }]0,\ov T[\x \ov  \Omega .
        \end{cases}
\end{equation}
It follows from \eqref{dimensionless_watertable} that  for $\ov t\in]0,\ov T[$ and $(\ov x,\ov z)\in\Omega^-_{\hab_0}(\ov t)$
\begin{equation}\label{pressure_general1}
 \disp\ov P_0(\ov t,\ov x,\ov z)=\rho\,g\big( \ov H_0(\ov t,\ov x)-\ov z \big), \quad 
 \disp\ov P_k(\ov t,\ov x,\ov z)=\rho\,g \ov H_k(\ov t,\ov x) \    \forall k>0 .
\end{equation}
 Equation \eqref{dimensionless_interface} gives
\begin{equation}\label{h_general1}
 \disp \hbotb(\ov x) \leq  \ov\ha_0(\ov t,\ov x)\le\max\Big\{\min\Big\{ \ov H_0(\ov t,\ov x)-\frac{\Pbubb}{\rho\,g},\hmaxb(\ov x) \Big\},\hbotb(\ov x) \Big\}\qquad\text{in }]0,\ov T[\x \ov \Omega_x .
\end{equation}
For the boundary conditions, we infer from the second and third equations of \eqref{dimensionless_transition} and from the second equation of \eqref{dimensionless_hydraulic} that,  for all $k\in\N$, 
\begin{equation}\label{boundary_general1}
 \begin{cases}
 \disp \alpha\,\ov P_k+\beta\,\ov u_k\cdot e_3=\ov F_k&\text{on } \domtb\x \Gammatb ,\\
  \disp\ov  P_0\big(\ov t,\ov x,\ov\ha_0(\ov t,\ov x)\big)=\rho\,g\big(\ov H_0(\ov t,\ov x)-\ov\ha_0(\ov t,\ov x) \big)&\text{for } \ov t\in]0,\ov T[,\quad \ov x\in \interfaceb(\ov t) ,\\
 \disp \ov K(\ov H_0) \nabla_{\ov x}\ov H_0 \cdot\ov  n=0&\text{on } \domtb\x \Gammavb .
\end{cases}  
\end{equation}
By \eqref{pressure_general1} for $k=1$,  $\partial_z \ov P_1=0$ on $\Omega^-_{\hab_0}(\ov t)$. Then by \eqref{velocity_general1} and the first equation of \eqref{pressure_general1}
\begin{equation}\label{u1}
\ov u_1 = 0\qquad\text{in }\Omega^-_{\hab_0}(\ov t).
\end{equation}

\paragraph{Short time case.} In this part,  $T=\ov T$, that is $\gamma=0$. The first equations of \eqref{dimensionless_transition} and \eqref{dimensionless_hydraulic}  become
  \begin{equation}\label{short_rescaled_coupled}
\begin{cases}
\disp \phi \deriv{s(\ov\Pplus)}{t}+\deriv{}{z}\big(\ov\vplus\cdot e_3\big)=0\phantom{\int_\int} & \text{for \ }\ov t\in]0,\ov T[\ ,\quad (\ov x,\ov z)\in\Omplusb(\ov t) ,\\
  \eps^2\dv_x\big(\ov K(\ov H)\,\nabla \ov H\big) = (\ov u_0\cdot e_3)|_{\interfaceb} &\text{for \ }(\ov t,\ov x)\in]0,\ov T[\x\ov \Omega_x .
\end{cases}
\end{equation}
We identify the main order terms appearing when the asymptotics \eqref{asymptotic_richards}--\eqref{asymptotic_sat} are substituted in the previous equations:
 for $\ov t\in]0,\ov T[$ and $(\ov x,\ov z)\in\Omega^+_{\ov\ha_0}(\ov t)$ 
\begin{equation}\label{order_0_coupled_1}
\phi \deriv{s(\ov P_0)}{t}+\deriv{}{ \ov z}\big(\ov u_0\cdot e_3\big)=0,
\end{equation}
\begin{equation}\label{order_0_coupled_2}
(\ov u_0 \cdot e_3)|_{\Gamma_{\ov\ha_0}}=0\qquad\text{for \ }(\ov t,\ov x)\in]0,\ov T[\x\ov \Omega_x.
\end{equation}
From \eqref{velocity_global} and \eqref{pressure_general1} we also compute $\ov u_0=0$  in $\Omega^-_{\ov\ha_0}(\ov t)$. In addition, from \eqref{h_general1} we get $s(\ov P_0)=1$ in $\Ommoinsb(\ov t)$ so that  $(\ov P_0,\ov u_0)$ satisfies \eqref{order_0_coupled_1} also in $\Ommoinsb(\ov t)$. The continuity of  $\ov u_0\cdot e_3$ being ensured by \eqref{order_0_coupled_2}, 
 $(\ov P_0,\ov u_0)$ satisfies \eqref{order_0_coupled_1}  in the whole $\Omega$.
By using \eqref{velocity_global} and \eqref{boundary_general1} we obtain the system \eqref{effective_short} and then the first claim of Proposition \ref{prop_richards} holds once again.

 \paragraph{Intermediate time case.}In this part,  $T=\eps^{-1}\ov T$, $\gamma=1$. 
 The first equation of \eqref{dimensionless_transition} and the equation \eqref{dimensionless_hydraulic_alternate}  become
  \begin{equation}\label{inter_rescaled_coupled}
\begin{cases}
\disp \phi \eps\deriv{s(\ov\Pplus)}{t}+\deriv{}{z}\big(\ov\vplus\cdot e_3\big)=0\phantom{\int_\int} & \text{for \ }\ov t\in]0,\ov T[\ ,\quad (\ov x,\ov z)\in\Omplusb(\ov t)\\
    \begin{aligned}
 &-\eps^2\dv_x\big(\ov K(\ov H)\,\nabla_x \ov H\big) 
 =-(\ov u\cdot e_3)|_{\Gammatb}-\eps\deriv{}{t}\Big(\int_{\hbotb(t,x)}^{\hsoilb(x)} \phi\,s(\ov\Pplus) \,dz \Big)
\end{aligned}&\text{for \ }(\ov t,\ov x)\in]0,\ov T[\x\ov \Omega_x
\end{cases}
\end{equation}
The corresponding main order relations are
  \begin{equation}\label{ue3_soil}
\ov u_0\cdot e_3=0\qquad \text{on } ]0,\ov T[\x\Gammatb
\end{equation}
  and for $\ov t\in]0,\ov T[$ and $ (\ov x,\ov z)\in\Omplusb(\ov t)$,
\begin{equation}\label{inter_order_0_coupled_1}
\deriv{}{z}\big(\ov u_0\cdot e_3\big)=0 .
\end{equation}
It follows that the constant vertical component of the velocity $\ov u_0\cdot e_3$ equals zero in $\Omplusb(\ov t)$.
We deduce from the first equation of \eqref{velocity_general1} that the pressure $\ov P_0$ is affine with respect to the $z$ variable with the slope $-\rho\,g$ in  $\Omplusb(\ov t)$. Accordingly, thanks to the first equation of \eqref{pressure_general1} and the continuity condition given in \eqref{boundary_general1},  the first equation of \eqref{effective_instant} holds. Using  relation  \eqref{velocity_global} we obtain the second equation of  \eqref{effective_instant}. 
Next, thanks to $\ov u_0=0$ and to the first equation of \eqref{boundary_general1} for $k=0$, we get 
 $\alpha\,\ov P_0=\ov F_0$.

If $\alpha\neq0$ then for all $(t,x)\in \domt\x\ov \Omega_x$ we have $P_0(\ov t,\ov x,\hsoilb(\ov t,\ov x))=\ov F_0(\ov t,\ov x)/\alpha$. Accordingly, thanks to the first equation of  \eqref{effective_instant}, we have
 \begin{equation*}
\ov H_0(\ov t,\ov x)=\frac{\ov F_0(\ov t,\ov x)}{\alpha\,\rho\,g}+\hsoilb(\ov t,\ov x).
\end{equation*}
The second claim of Proposition \ref{prop_richards} in the case $\alpha\neq0$ is proved.

If $\alpha=0$, the compatibility condition $\ov F_0=0$ is imposed.
 After identifying the coefficients associated with $\eps^1$  in the second  equation of  \eqref{inter_rescaled_coupled} we have
\begin{equation*}
0=-(\ov u_1\cdot e_3)|_{\Gammatb}-\deriv{}{t}\Big(\int_{\hbotb(x)}^{\hsoilb(x)} \phi\,s(\ov P_0) \,d\ov z \Big) 
\end{equation*}
and, with the first equation of \eqref{effective_instant}, 
\begin{equation*}
\rho\,g\,\Big(\int_{\hbotb(x)}^{\hsoilb(x)} \phi\,s'(\ov P_0) \,d\ov z \Big)\deriv{\ov H_0}{t}=-( \ov u_1\cdot e_3)|_{\Gammatb} .
\end{equation*}
The first equation of \eqref{boundary_general1} for $k=1$ implies, since $\alpha=0$, that $(\ov u_1\cdot e_3)|_{\Gammatb}=\ov F_1/\beta$. This ends the proof of the second claim of Proposition \ref{prop_richards}  in the case $\alpha=0$.

 \paragraph{Long time case.}In this part, $T=\eps^{-\gamma}\ov T$, $\gamma=2$.
The first equation of  \eqref{dimensionless_transition} and equation \eqref{dimensionless_hydraulic_alternate} are now
   \begin{equation}\label{long_rescaled_coupled}
\begin{cases}
\disp \phi\, \eps^2\deriv{s(\ov\Pplus)}{t}+\deriv{}{z}\big(\ov\vplus\cdot e_3\big)=0 & \text{for \ }\ov t\in]0,\ov T[\ ,\quad (\ov x,\ov z)\in\Omplusb(\ov t)\\
 \begin{aligned}
 &- \eps^2\,\dv_{\ov x}\Big(\ov K(\ov  H)\,\nabla_{\ov x}\ov H \Big) 
 =-(\ov u\cdot e_3)|_{\Gammatb}- \eps^2\deriv{}{\ov t}\Big(\int_{\hbotb(\ov x)}^{\ov \hsoil(\ov x)} \phi\,s(\ov \Pplus) \,d\ov z \Big)
\end{aligned}&\text{for \ }(\ov t,\ov x)\in]0,\ov T[\x\ov \Omega_x \\
\end{cases}
\end{equation}
As in the intermediate time case, we substitute asymptotics \eqref{asymptotic_richards}--\eqref{asymptotic_sat} in the previous equations. 
Identifying the coefficients associated with $\eps^n$ for $n\in\{1,2\}$, we get $\partial_z(\ov u_n\cdot e_3)=0$ in $\Omplusb(\ov t)$ and $\ov u_n\cdot e_3=0$ on $\Gammatb$. This leads to  
\begin{equation}\label{u_10}
\ov u_0\cdot e_3=\ov u_1\cdot e_3=0\qquad\text{on }\Omplusb(\ov t).
\end{equation} By using the same arguments we obtain $\ov P_0 = \rho\,g(\ov H_0-z)$ and $\ov v_0=0$ in whole $\Omega$.  System \eqref{effective_instant} is satisfied. 
The characterization of $\ov H_0$ depends on the values of $\alpha$. Similar arguments to those employed  in the intermediate time case when $\alpha\neq0$ lead to \eqref{effective_alpha_0}.

It remains to deal with the case $\alpha=0$. 
In this case we first remark that the compatibility condition $\ov F_0=0$ holds (see
\eqref{boundary_general1} for $k=0$). Furthermore,  since $\ov P_0=\rho\,g(\ov H_0+z)$ we get from \eqref{u1} and \eqref{u_10} that $\ov u_1=0$ in $]0,\ov T[\x \ov \Omega$. 
Thus, using
\eqref{velocity_general2} and \eqref{velocity_global} we get $\ov v_1=\ov w_1=-k_r(\ov P_0) \ov M_0\,\nabla_x \ov H_0$. Moreover the first equation of \eqref{boundary_general1} for $k=1$ gives $\ov F_1=0$ (since $\alpha=0$).
It remains to get the first relation of system  \eqref{effective_long}. By plugging asymptotics \eqref{asymptotic_richards}--\eqref{asymptotic_sat} in  the second equation of \eqref{long_rescaled_coupled} and by identifying the coefficients associated with $\eps^2$ we get (\eqref{velocity_general1})
    \begin{equation}\label{long_order_2_coupled}
 -\dv_x\big( \ov K(\ov H_0)\,\nabla_x\ov H_0 \big) =- (\ov u_2\ cdot e_3)|_{\Gammatb}-\deriv{}{t}\Big(\int_{ \hbotb(x)}^{\hsoilb(x)} \phi\,s(\ov P_0) \,dz \Big)
\qquad\text{for \ }(\ov t,\ov x)\in]0,\ov T[\x\ov \Omega_x .
\end{equation}
We end the proof by noting that, thanks to the equality $\alpha=0$ and the first equation of \eqref{boundary_general1} for $k=2$, we have $( \ov u_2\cdot e_3)|_{\Gammatb}=\ov F_2/\beta$.
\qed

\end{document}